\documentclass[12pt]{article}

\usepackage{geometry}
\usepackage[T1]{fontenc}
\usepackage{lmodern}
\usepackage{microtype}
\usepackage{mathtools,amssymb,amsthm,mathrsfs,bbm}
\usepackage{enumitem}
\usepackage{xcolor}
\usepackage{booktabs}
\usepackage{hyperref}
\usepackage[nameinlink,noabbrev]{cleveref}

\hypersetup{
  colorlinks=true,
  linkcolor=blue!55!black,
  citecolor=green!40!black,
  urlcolor=blue!60!black,
  pdftitle={Todd's relation conjecture and binary relations for multiple zeta values in positive characteristic},
  pdfsubject={Multiple zeta values and Carlitz multiple polylogarithms in positive characteristic}
}

\setlist[enumerate,1]{label=\textup{(\roman*)},leftmargin=2.2em}
\setlist[itemize]{leftmargin=1.8em}
\numberwithin{equation}{section}

\newtheorem{theorem}{Theorem}[section]

\newtheorem{lemma}[theorem]{Lemma}
\newtheorem{corollary}[theorem]{Corollary}

\theoremstyle{definition}
\newtheorem{definition}[theorem]{Definition}

\newtheorem{remark}[theorem]{Remark}
\newtheorem*{acknowledgements}{\bf Acknowledgements}

\allowdisplaybreaks[1]

\DeclareMathOperator{\Span}{Span}

\newcommand{\Fq}{\mathbb F_q}

\newcommand{\Li}{\mathrm{Li}}
\newcommand{\Si}{\mathrm{Si}}
\newcommand{\rS}{\mathrm{S}}

\newcommand{\ind}{\operatorname{Ind}}
\newcommand{\BCP}{\mathcal{BCP}}
\newcommand{\cB}{\mathcal B}
\newcommand{\cP}{\mathcal P}
\newcommand{\cBs}{\mathcal B^{\ast}}
\newcommand{\cC}{\mathcal C}
\newcommand{\sC}{\mathscr C}
\newcommand{\cL}{\mathcal L}
\newcommand{\cH}{\mathcal H}
\newcommand{\sL}{\mathscr L}

\newcommand{\cR}{\mathscr R}
\newcommand{\cS}{\mathscr S}

\newcommand{\cZ}{\mathcal Z}
\newcommand{\cI}{\mathcal I}
\newcommand{\cE}{\mathcal E}
\newcommand{\br}{\mathfrak{BR}}

\newcommand{\bbZ}{\mathbb Z}

\newcommand{\wt}{\operatorname{wt}}
\newcommand{\dep}{\operatorname{dep}}
\newcommand{\fs}{\mathfrak{s}}
\newcommand{\fh}{\mathfrak{h}}

\newcommand{\fm}{\mathfrak{m}}
\newcommand{\ft}{\mathfrak{t}}
\newcommand{\fu}{\mathfrak{u}}
\newcommand{\fv}{\mathfrak{v}}

\newcommand{\fa}{\mathfrak{a}}
\newcommand{\fb}{\mathfrak{b}}

\newcommand{\fx}{\mathfrak{x}}
\newcommand{\fy}{\mathfrak{y}}

\newcommand{\BTTw}{\mathcal{B}_w^{(\textup{TT})}}
\newcommand{\Ss}{\ast}
\newcommand{\St}{\triangleleft}
\newcommand{\Sc}{\diamond}
\newcommand{\bQ}{\boldsymbol{Q}}
\newcommand{\cD}{\mathscr D}
\newcommand{\im}{\operatorname{Im}}

\title{On the $\Fq$-linear relations among multiple zeta values in positive characteristic}
\author{Jinyuan Hu}
\date{}

\begin{document}

\maketitle

\begin{abstract}
This paper gives a new path to determine all $\Fq$-linear relations among multiple zeta values in positive characteristic, which had been done by Huang, Lai, Li and the author in \cite{HL2026}. The path used in this paper is based on the idea of operator generation, which originated from the work of Todd in \cite{Todd2018}, developed by the author, and independently by Im–Kim–Ngo Dac in \cite{IKN2026}.
\end{abstract}

\tableofcontents

\section{Introduction}
For any non-empty set $A$, we denote by $A^{\bullet}$ the disjoint union $\{\emptyset\} \sqcup_{k=1}^{\infty} A^k$. 
Let $\bbZ_+$ be the set of positive integers, and denote by $\cI$ the set $\bbZ_+^{\bullet}$.
The elements in $\cI$ are called \emph{tuples}.
For any tuple $\fs =(s_1,\cdots,s_r) \in \cI$, the \emph{depth} $\dep(\fs)$ and \emph{weight} $\wt(\fs)$ are defined by $\dep(\fs) = r$ and $\wt(\fs)=\sum_{i=1}^{r} s_i$, with the convention $\dep(\emptyset)=0$ and $\wt(\emptyset)=0$. For a tuple $\fs=(s_1,\cdots,s_l)$, we set $\fs=(s_1,\fs^T)=(\fs^H,s_l)$ (possibly $\fs^H,\fs^T=\emptyset$). For any $n\in\bbZ_+$, we set $\prescript{n+}{}{\fs}=(s_1+n,\fs^T)$ and $\fs^{+n}=(\fs^H,s_l+n)$

Define $\cI_w = \left\{ \fs \in \cI \mid \wt(\fs) = w \right\}$ for any weight $w \ge 0$, and define $\cI_{>0} = \cI \setminus \{\emptyset\}$.
For $w\in\bbZ_+$ and a field $F$, we denote the $F$-linear formal sum space of tuples of weight $w$ by $\cS_w^{(F)}=\bigoplus_{\fs\in\cI_w} F\fs$, and $\cS^{(F)}=\bigoplus_{\fs\in\cI} F\fs$.
As usual, $[n]$ denotes the set $\{1,2,\ldots,n\}$ for any $n \in \bbZ_+$, and $[0]=\emptyset$.

\subsection{Thakur's multiple zeta values and Carlitz multiple polylogarithm values}

Let $A=\Fq[\theta]$ be the polynomial ring in one variable $\theta$ over a finite field $\Fq$ of characteristic $p$, and let $A_+$ be the set of monic polynomials in $A$. Let $K=\Fq(\theta)$ be its fraction field, and $K_\infty=\Fq((\theta^{-1}))$ be the completion of $K$ at $\infty$.

For a tuple $\fs = (s_1,\ldots,s_r) \in \cI_{>0}$ and $d\in\bbZ_{\geq0}$, define the power sum $\rS_d(\fs)$ as
\[
 \rS_d(\fs)=
 \sum_{\substack{a_1,\ldots,a_r\in A_+\\
                  d=\deg a_1>\cdots>\deg a_r}}
 \frac{1}{a_1^{s_1}\cdots a_r^{s_r}}\in K.
\]
For convenience, we set $\rS_{-1}(\emptyset)=1$, $\rS_{d}(\emptyset)=0$ for $d\in\bbZ\setminus\{-1\}$, and $\rS_d(\fs)=0$ for $d<\dep(\fs)-1$.

We further define
\[
\rS_{<d}(\fs)=\sum_{k<d}\rS_k(\fs).
\]
 then
\[
\rS_d(\fs)=\rS_d(s_1)\rS_{<d}(\fs^T).
\]

Now we define the multiple zeta values in positive characteristic (abbreviated as MZV) introduced by Thakur \cite{Thakur2004} as follows:
\[
\zeta(\fs)=\sum_{d\geq0}\rS_d(\fs)\in K_{\infty}.
\]

For $\rS_d(\fs), \rS_{<d}(\fs)$ and $\zeta(\fs)$, we call $\dep(\fs)$ and $\wt(\fs)$ the depth and weight of them.

The other object closely associated with MZVs is Carlitz multiple polylogarithm values. Set $l_0=1$ and $l_d=\prod_{k=1}^d(\theta-\theta^{q^k})$ for $d\in\bbZ_+$, for a tuple $\fs=(s_1,\cdots,s_r)\in\cI_w$, we introduce analogues of power sums, which we called Carlitz power sum here
\[
\Si_d(\fs)=\sum_{d=k_1>\cdots>k_r\geq0}\frac{1}{l_{k_1}^{s_1}\cdots{l_{k_r}^{s_r}}}\in K_{\infty},
\]
Similarly, we set $\Si_{-1}(\emptyset)=1$, $\Si_{d}(\emptyset)=0$ for $d\in\bbZ\setminus\{-1\}$, and $\Si_d(\fs)=0$ for $d<\dep(\fs)-1$.

Further define
\[
\Si_{<d}(\fs)=\sum_{k<d}\Si_k(\fs)\in K_{\infty}.
\]
Similarly, we have
\[
\Si_d(\fs)=\Si_d(s_1)\Si_{<d}(\fs^T).
\]

The Carlitz multiple polylogarithm values (abbreviated as CMPLVs) are defined by
\[
\Li(\fs)=\sum_{d\geq0}\Si_d(\fs)\in K_{\infty}.
\]
For $\Si_d(\fs), \Si_{<d}(\fs)$ and $\Li(\fs)$, we still call $\dep(\fs)$ and $\wt(\fs)$ the depth and weight of them.

For MZVs, we are interested in their linear relations, over $K$ or over $\Fq$. However, the linear relations among CMPLVs are the objects which are easier to treat in general. \cite{Hu2026} explains that studying the linear relations among MZVs is equivalent to studying the ones among CMPLVs. More precisely, we have the following theorem.

\begin{theorem}[Hu \cite{Hu2026}]\label{thm_MZVCMPLV}
    Define $\zeta$ and $\Li$ as the $F$-linear maps from $\cS^{(F)}$ to $K_{\infty}$ by linear extension for $F=K$ or $\Fq$. Then there is a canonical linear automorphism $T$ of $\cS^{(F)}$, such that $\zeta\circ T=\Li$.
\end{theorem}
Here ``canonical'' means $T$ preserves some additional algebraic structures, but we will not use this property in this paper.

Note that the kernels of $\zeta$ and $\Li$ are actually the $F$-linear relation space of MZVs and CMPLVs, so \Cref{thm_MZVCMPLV} enables us to focus on the linear relations among CMPLVs.

For the $K$-linear relations, Todd and Thakur predicted an explicit basis and dimension for the
$K$-span of CMPLVs, which is now a theorem, independently due
to Chang--Chen--Mishiba \cite{CCM2023} and Im--Kim--Le--Ngo Dac--Pham \cite{IKLNDP2024}. 

Define $\cZ_{w}^{(K)}$ the $K$-linear subspace of $K_\infty$ spanned by CMPLVs $\Li(\mathfrak{s})$ of weight $w$.

\begin{theorem}[Chang--Chen--Mishiba \cite{CCM2023}, independently Im--Kim--Le--Ngo Dac--Pham \cite{IKLNDP2024}]\label{thm_CCMIKLNP}
For any weight $w \ge 1$, we have $\dim_{K} \cZ_w^{(K)} = d_w^{(K)}$, where the sequence $\{d_w^{(K)}\}_{w \ge 1}$ is defined by the generating function 
\[
\sum_{w=1}^{\infty} d_w^{(K)} x^w = \frac{x(1-x^{q-1})}{1-2x+x^{q+1}}.
\]
Moreover, $\left\{ \Li(\fs) ~|~ \fs \in \BTTw  \right\}$ is a $K$-basis of $\cZ_w^{(K)}$, where 
\[
\BTTw = \left\{ (\fh,z) \in \cI_w ~|~ \fh \in [q]^{\bullet}, z \in [q-1] \right\}.
\]
\end{theorem}

Todd also formulated a conjecture about how to obtain $K$-linear relations, which is the idea the author called ``operator generation''. He claimed all of $K$-linear relations are generated by lifting the binary relations from applying the $\cB,\cBs$ and $\cC$ operators to the fundamental relation, which we will introduce later. This conjecture is also a theorem now, due to Chang–Chen–Mishiba \cite{CCM2023} and the author \cite{Hu2026}.

We first define binary relations.

\begin{definition}
    \begin{enumerate}
        \item A binary relation $R$ of CMPLVs of weight $w$ is given by a collection of elements $a_i,b_j$ of $K$ such that for all $d\in\bbZ$
        \[
        \sum_ia_i\Si_d(\fs_i)+\sum_jb_j\Si_{d+1}(\ft_j)=0,
        \]
        where the sum runs through tuples $\fs_i$ and $\ft_j$ of weight $w$. We denote the above equality by $R(d)$. And we may denote $R$ by $[P;Q]$ where $P=\sum_ia_i\fs_i$ and $Q=\sum_jb_j\ft_j$.
        \item A binary relation is called a fixed relation if $b_j=0$ for all $j$.
    \end{enumerate}
\end{definition}

We denote the set of all binary relations of $\Si$ of weight $w$ by $\br_w$.

The simplest binary relation, in the sense of owning the lowest weight, is the one we called the fundamental relation, denoted by $R_1$ (see \cite{Thakur2009}, §3.4.6):
\[
    \Si_d(q)+D_1\Si_{d+1}(1,q-1)=0,\quad \text{where } D_1=\theta^q-\theta.
\]

Then we introduce the $\cB,\cBs$ and $\cC$ operators.

\begin{definition}
    Given a tuple $W=(w_1,\cdots,w_r)\in\cI_k$.
    \begin{enumerate}
        \item The operator $\cB_W:\br_w\to\br_{w+k}$ is a linear map defined by
        \[
            \cB_W(R): \quad\cB_W(R)(d)=\Si_d(W)\sum_{j<d}R(j) \quad \text{for any } R\in\br_w.
        \]
        \item The operator $\cBs_W:\br_w\to\br_{w+k}$ is a linear map defined by 
        \[
        \cBs_W=\cB_{(w_1)}\circ\cdots\circ\cB_{(w_r)}   . 
        \]
        \item The operator $\cC_W:\br_w\to\br_{w+k}$ is a linear map defined by
        \[
            \cC_W(R): \quad\cC_W(R)(d)=R(d)\Si_{<d+1}(W) \quad \text{for any } R\in\br_w.
        \]
    \end{enumerate}
    For $W=\emptyset$, we set all these operators to be the identity.
\end{definition}

Defining the product maps $\Ss,\St$ and $\Sc$ would help us understand $\cB,\cBs$ and $\cC$ better.
    \begin{align}
        \fs\Ss\ft=&(s_1,\fs^T\Ss\ft)+(t_1,\fs\Ss\ft^T)+(s_1+t_1,\fs^T\Ss\ft^T);\label{eq:Ss}\\
        \fs\St\ft=&(s_1,\fs^T\Ss\ft);\\
        \fs\Sc\ft=&(s_1+t_1,\fs^T\Ss\ft^T).
    \end{align}
    where $\fs\Ss\emptyset=\emptyset\Ss\fs=\fs$, $\fs\Sc\emptyset=\emptyset\Sc\fs=0t$ and $\fs\St\emptyset=\fs$. And $\Ss$ and $\Sc$ are both commutative and associative. For multiple $\St$-products, we always evaluate from right to left.

\begin{theorem} For any $\fs,\ft\in\cI$, and $d\in\bbZ_{\geq0}$,
    \begin{align*}
        \Si_{<d}(\fs)\Si_{<d}(\ft)&=\Si_{<d}(\fs\Ss\ft);\\
        \Si_d(\fs)\Si_d(\ft)&=\Si_d(\fs\Sc\ft);\\
        \Si_d(\fs)\Si_{<d}(\ft)&=\Si_{d}(\fs\St\ft),
    \end{align*}
\end{theorem}
\begin{proof}
    Verify directly by definition
\end{proof}

Then the $\cB,\cBs,\cC$ can be rewritten as
\begin{align*}
    \cB_W[P; Q]&=[W\St P + W\Sc Q+W\St Q; 0];\\
    \cBs_W[P; Q]&=\Bigl[(w_1)\St(w_2)\St\cdots\St\bigl((w_r)\St P + (w_r)\Sc Q+(w_r)\St Q\bigr); 0\Bigr];\\
    \cC_W[P; Q]&=[P\Sc W+P\St W; Q\St W].
\end{align*}

Finally, we need to explain what is the lift of a binary relation. Given a binary relation $R=[P;Q]$, it is
\[\Si_d(P)+\Si_{d+1}(Q)=0,\]
summing over all d, we obtain a relation among CMPLVs,
\begin{equation}\label{eq:lift}
    \Li(P)+\Li(Q)=0.
\end{equation}

We denote $P+Q$ by $\overline{R}$ or $\overline{[P;Q]}$, call it the lift of $R$. Then \Cref{eq:lift} is actually $\Li(\overline{R})=0$.

Now we can state Todd's relation conjecture.
\begin{theorem}[Chang–Chen–Mishiba \cite{CCM2023}, and Hu \cite{Hu2026}]\label{thm:Toddconj}
    For $n\in\bbZ_{+}$, the $K$-linear relation space among CMPLVs of weight $q+n$ is
spanned by
\begin{equation}\label{eq:Toddconj}
 \{\overline{\cB_U(R_1)}: U\in\cI_n\}\cup\{\overline{\cC_V(R_1)}: V\in\cI_n\}\cup
 \{\overline{\cB_U\circ\cC_V(R_1)}:U\in\cI_l, V\in\cI_m, l+m=n\}.
\end{equation}
It is also spanned when every $\cB$ in \eqref{eq:Toddconj} is replaced
by $\cBs$.
\end{theorem}
\begin{remark}
    \Cref{thm:Toddconj} was originally formulated for MZVs. The $\cB$-parts of this theorem are equivalent for MZVs and CMPLVs. But $\cBs$-parts need to be treated separately. See \cite{Hu2026} for more details.
\end{remark}

\subsection{Main results in this paper}
Now we turn to the $\Fq$-linear relation among CMPLVs. Huang, Lai, Li and the author have obtained a complete description in \cite{HL2026}. Before introducing the detailed results, we introduce the quadruple-carry relations constructed by them.

\begin{lemma}
    For any $\fh,\fm,\ft \in \cI$, $r_1 \in \bbZ_+$, and $r_2 \in \bbZ_{\geqslant 0}$ satisfying either $r_2>0$, or $r_2=0$ and $\ft=\emptyset$ (with the convention that, in the latter case, $(r_2,\ft)=(0,\emptyset)$ is identified with $\emptyset$), define the following formal sum of tuples:
\begin{align*}
& \bQ\left( \fh,\boxed{q+r_1},\fm,\boxed{q+r_2},\ft \right) \notag\\
\coloneq & \bigl( \fh, (q+r_1,\fm)^{+1}, (q-1)*(r_2,\ft) \bigr) + \bigl( \fh, q+r_1, \fm, 1, (q-1)*(r_2,\ft) \bigr) \notag\\
&+ \bigl( \fh, q, (r_1,\fm)^{+1}, (q-1)*(r_2,\ft) \bigr) + \bigl( \fh, q, r_1, \fm, 1, (q-1)*(r_2,\ft) \bigr) \notag\\
&- \mathbbm{1}_{\fh \neq \emptyset} \cdot \bigl( \fh^{+1},(q-1)*(r_1,\fm,q+r_2,\ft) \bigr) - \bigl( \fh,1, (q-1)*(r_1,\fm,q+r_2,\ft) \bigr) \notag\\
&- \mathbbm{1}_{\fh \neq \emptyset, r_2 \neq 0} \cdot \bigl( \fh^{+1},(q-1)*(r_1,\fm,q,r_2,\ft) \bigr) - \mathbbm{1}_{r_2 \neq 0} \cdot \bigl( \fh,1,(q-1)*(r_1,\fm,q,r_2,\ft) \bigr). \label{def_Q}
\end{align*}
Then
\[\label{eq:qcrelation}
    \Li\bigl(\bQ\left( \fh,\boxed{q+r_1},\fm,\boxed{q+r_2},\ft \right)\bigr) = 0.
\]
We call the above equation a quadruple-carry relation.
\end{lemma}
\begin{proof}
    See \cite{HL2026}, Lemma 2.5.
\end{proof}

The following is the main result in \cite{HL2026}.

\begin{theorem}[Hu--Huang--Lai--Li \cite{HL2026}]\label{thm:HLmain}
 Define $\cZ_{w}$ the $\Fq$-linear subspace of $K_\infty$ spanned by CMPLVs of weight $w$, and $\cR_w$ the $\Fq$-linear subspace of $\cS_w^{(\Fq)}$ spanned by the $\Fq$-linear relations among CMPLVs of weight $w$ (i.e. the intersection of $\operatorname{Ker}\Li$ and $\cS_w^{(\Fq)}$). For any weight $w\geq1$,
    \begin{enumerate}
        \item $\dim_{\Fq} \cZ_w = d_w$, where the sequence $\{d_w\}_{w \geqslant 1}$ is defined by the generating function
        \[\label{eq:dwgf}
            \sum_{w=1}^{\infty} d_w x^w = \frac{x(1-x^q)(1-2x+x^q)}{(1-2x+x^{q+1})^2}.
        \]
        \item $\left\{ \Li^{\star}(\fs) ~\big|~ \fs \in \cB_w  \right\}$ is an $\Fq$-basis of $\cZ_w$, where $\cB_w = \cB \cap \cI_w$,
        \[\label{eq:hlbasis}
            \cB = [q]^{\bullet} \bigcup \left\{ (\fh,m,\ft) \in \cI ~\Big|~ \fh,\ft \in [q]^{\bullet}, m \in [2q] \setminus [q] \right\},
        \]
        and \[\Li^{\star}(\fs) = \sum_{\substack{n_1\geqslant n_2 \geqslant \cdots \geqslant n_k \\ n_1,n_2,\ldots,n_k \in \bbZ_{\geqslant 0}}}  \frac{ 1 }{\ell_{n_1}^{s_1} \ell_{n_2}^{s_2} \cdots \ell_{n_k}^{s_k}}.\]
        \item $\cR_w$ is the $\Fq$-span of all quadruple-carry relations of weight $w$.
    \end{enumerate}
\end{theorem}

The key point of \cite{HL2026} is actually \Cref{thm:HLmain}~(iii), which completely describes $\Fq$-linear relation space $\cR_w$. However, this is not the unique way to understand $\cR_w$. 

While developing the work \cite{HL2026}, inspired by \Cref{thm:Toddconj}, the author wanted to use the similar operator generation method to study $\cR_w$. Later, Im, Kim and Ngo Dac formulated a similar idea in \cite{IKN2026}. Using the terminology of \cite{IKN2026}, we first introduce the so-called $\cP$ operator.

For any $\fs=(s_1,\cdots,s_l)\in\cI\setminus\{\emptyset\}$, define the $\Fq$-linear operator $\cP:\cS_w^{(\Fq)}\to\cS_{w+2q}^{(\Fq)}\oplus\cS_{w+2q}^{(\Fq)}$ by
\begin{align*}
     \cP(\fs)=\biggl[(q)(\St+\Sc)\Bigl(\bigl(\fs^{+1}+(\fs,1),q-1\bigr)\Bigr);-\bigl(1,(q-1)\Ss(\fs,q)\bigr)\biggr],
\end{align*}
where $\fs(\St+\Sc)\ft=\fs\St\ft+\fs\Sc\ft$.

The first result in this paper shows that $\cBs,\cC, \cP$ operators can generate all $\Fq$-linear relations. We will prove it in two different ways in \Cref{sec:proof} and \Cref{sec:another} respectively.

\begin{theorem}\label{thm:BCP}
    $\cR_w=\Span_{\Fq}\{\overline{\cBs_{\fu}\circ\cC_{\fv}\circ\cP(\fs)} ~|~ \wt(\fu)+\wt(\fv)+\wt(\fs)+2q=w\}$
\end{theorem}

We will use $\BCP_w$ to denote the RHS of the above equation.

\begin{remark}
    Similarly, we can use $\cB$ to replace $\cBs$. See \cite{Hu2026}, Theorem 3.3.
\end{remark}

We will also calculate $\dim_{\Fq}\cR_w$ based on \Cref{thm:BCP} in \Cref{sec:dimension}.
\begin{theorem}\label{thm:dimension}
    Let $\dim_{\Fq}\cR_w=r_w$, we have
    \[\label{rwgf}
    \sum_{w=1}^{\infty} r_w x^w = \frac{x^{2q+1}(1-x)^2}{(1-2x)(1-2x+x^{q+1})^2}.
\]
\end{theorem}

\subsection{List of notations}
Some of notations in this paper are unnecessary in fact, but we still use them for convenience. Here we list main notations that may be used repeatedly in different subsections.
\begin{align*}
    |A|&=\text{the cardinality of a set }A;\\
    A^{\bullet}&=\{\emptyset\} \sqcup_{k=1}^{\infty} A^k;\\[6pt]
    \cI&=\bbZ_+^{\bullet};\\
    [n]&=\{1,2,\cdots,n\};\\
    K&=\Fq(\theta)=\text{the rational function field over the finite field }\Fq;\\
    D_1&=\theta^q-\theta;\\
    K_\infty&=\Fq((\theta^{-1}))=\text{the completion of }K\text{ at infinite place};\\[6pt]
    \text{For a }\fs=(s_1,\cdots,s_l)\in\cI:\\
    \dep(\fs)&=l;\\
    \wt(\fs)&=\sum_{i=1}^ls_i;\\
    \prescript{n+}{}{\fs}&=(n+s_1,\cdots,s_l);\\
    \fs^{+n}&=(s_1,\cdots,s_l+n);\\
    \fs^H&=(s_1,\cdots,s_{l-1});\\
    \fs^T&=(s_2,\cdots,s_l);\\[6pt]
    \cI_w&=\{\fs\in\cI~|~\wt(\fs)=w\};\\
    \cS_w^{(F)}&=\bigoplus_{\fs\in\cI_w}F\fs\quad\text{for a field }F;\\
    \cS^{(F)}&=\bigoplus_{w\geq0}F\cS_w^{(F)}\quad\text{for a field }F;\\
    \rS_d(\fs)&=\text{the power sum of }\fs;\\
    \rS_{<d}(\fs)&=\sum_{k<d}\rS_k(\fs);\\
    \zeta{\fs}&=\text{the Thakur's multiple zeta value of }\fs;\\
    \Si_d(\fs)&=\text{the Carlitz power sum of }\fs;\\
    \Si_{<d}(\fs)&=\sum_{k<d}\Si_k(\fs);\\
    \Li(\fs)&=\text{the Carlitz multiple polylogarithm value of }\fs;\\
    \BTTw&= \left\{ (\fh,z) \in \cI_w ~|~ \fh \in [q]^{\bullet}, z \in [q-1] \right\};\\
    d_w^{(K)}&=|\BTTw|;\\
    \br_w&=\text{the set of binary relations with weight w}\\&=\{[P;Q]~|~P,Q\in\cS_w^{(K)},\Si_d(P)+\Si_{d+1}(Q)=0~\text{for all }d\};\\
    R_1&=\text{fundamental relation}=\bigl[(q);D_1(1,q-1)\bigr];\\[6pt]
    \text{The } \Ss,\St,\Sc \text{ products}:\\
    \fs\Ss\ft&=(s_1,\fs^T\Ss\ft)+(t_1,\fs\Ss\ft^T)+(s_1+t_1,\fs^T\Ss\ft^T);\\
    \fs\St\ft&=(s_1,\fs^T\Ss\ft);\\
    \fs\Sc\ft&=(s_1+t_1,\fs^T\Ss\ft^T);\\[6pt]
    \text{The } \cB,\cBs,\cC \text{ operators}:\\
     \cB_W[P; Q]&=[W\St P + W\Sc Q+W\St Q; 0];\\
    \cBs_W[P; Q]&=\Bigl[(w_1)\St(w_2)\St\cdots\St\bigl((w_r)\St P + (w_r)\Sc Q+(w_r)\St Q\bigr); 0\Bigr];\\
    \cC_W[P; Q]&=[P\Sc W+P\St W; Q\St W];\\[6pt]
    \overline{[P;Q]}&=P+Q\quad\text{the lift of }[P;Q];\\
    \cR_w&=\{P\in\cS_w^{(\Fq)}~|~\Li(P)=0\};\\
    \cP(\fs)&=\biggl[(q)(\St+\Sc)\Bigl(\bigl(\fs^{+1}+(\fs,1),q-1\bigr)\Bigr);-\bigl(1,(q-1)\Ss(\fs,q)\bigr)\biggr];\\
    \BCP_w&=\Span_{\Fq}\{\overline{\cBs_{\fu}\circ\cC_{\fv}\circ\cP(\fs)} ~|~ \wt(\fu)+\wt(\fv)+\wt(\fs)+2q=w\};\\
    r_w&=\dim_{\Fq}\cR_w;\\
    \cD_0,\cD_1&:~\cD_0,\cD_1\text{ operators, see \Cref{lem:D}};\\[6pt]
    \text{The }\cL,\cH\text{ operators}:\\
    \cL[\fs;\ft]&=(\fs,\prescript{q+}{}{\ft})+\mathbbm{1}_{\ft\neq\emptyset}(\fs,q,\ft);\\
    \cH[\fs;\ft]&=\mathbbm{1}_{\fs\neq\emptyset}\bigl(\fs^{+1},(q-1)\Ss\ft\bigr)+\bigl(\fs,1,(q-1)\Ss\ft\bigr);\\[6pt]
    \sL_w&=\Span_{\Fq}\{\cL[\fh;\fs] ~|~ \fh\in[q]^{\bullet},~ \wt(\fh)+\wt(\fs)+q=w\};\\
    \ind\bigl(\cL[\fs;\ft]\bigr)&=\cH[\fs;\ft];\\
    \fs^{i},\fs^{cano},E_i(\fs),E(\fs)&:~\text{see \Cref{subsec:KernelL}};\\
    \cP[\fs;\ft]&=\biggl[(q)(\St+\Sc)\bigl(\fs^{+1}+(\fs,1),(q-1)\Ss\ft\bigr);-\Bigl(1,(q-1)\Ss\bigl(\fs,(q)(\St+\Sc)\ft\bigr)\Bigr)\biggr];\\
    \sC&=\Span_{\Fq}\{[\fh;\ft] ~|~ \fh\neq\emptyset\};\\
    H_+&=\Span_{\Fq}\{\bigl(\fs,(a+1)+(a,1),(q-1)\Ss\ft\bigr) ~|~ \fs,\ft\in\cI,a\in\bbZ+\};\\
    \rho_n&:~\text{the extraction operator, see the proof of \Cref{lem:separate}};\\
    \Lambda[\fs;\ft]&=\bigl[(q-1)\Ss\fs;\ft\bigr]+\bigl[\fs;(q-1)\Ss\ft\bigr];\\
    U_w&=\{P+D_1Q~|~P,Q\in\cS_w^{(\Fq)}~,~\Li(P+D_1Q)=0\};\\
    V_w^{\text{Fix}}&=\{P\in V_w~|~[P;0]\in\br_w\}.
\end{align*}

\section{Completeness of the BCP relations}\label{sec:proof}
In this section, we prove \Cref{thm:BCP}. The core idea is very similar to the one of \cite{HL2026}, with detailed techniques differing a lot.
\subsection{The image of $\cP$ lies in $\br_{w+2q}$}
First step is proving that the image of $\cP$ lies in $\br_{w+2q}$. We first construct the $\cD_0$ and $\cD_1$ operator.

\begin{lemma}\label{lem:D}
    For any binary relation $R$ of the form $\Bigl[(q)(\St+\Sc)P;-\bigl(1,(q-1)\Ss Q\bigr)\Bigr]$, there are the following two binary relations, denoted by $\cD_0(R)$ and $\cD_1(R)$.
    \begin{align}
        &\cD_0(R)=\Bigl[(q)(\St+\Sc)(1,P);-\bigl(1,(q-1)\Ss (1,Q)\bigr)\Bigr];\label{eq:D0}\\
        &\cD_1(R)=\Bigl[(q)(\St+\Sc)\prescript{1+}{}{P};-\bigl(1,(q-1)\Ss \prescript{1+}{}{Q}\bigr)\Bigr];\label{eq:D1}
    \end{align}
    where $(1,P)=\sum_ia_i(1,\fs_i)$ and $\prescript{1+}{}{P}=\sum_ia_i\prescript{1+}{}{\fs_i}$ for any $P=\sum_i a_i\fs_i$.
\end{lemma}
\begin{proof}
    To prove \Cref{eq:D0}, we have
    \begin{align*}
        \Si_d(q)\Si_{<d+1}(1,P)&=\Si_d(q)\sum_{n<d+1}\Si_n(1,P)\\&=\Si_d(q)\sum_{n<d+1}\Si_n(1)\Si_{<n}(P)\\&=\Si_d(q)\sum_{n<d+1}\Si_n(1)\Si_{n-1}^{-1}(q)\Si_{n-1}\bigl((q)(\St+\Sc)P\bigr)\\&=\Si_d(q)\sum_{n<d+1}\Si_n(1)\Si_{n-1}^{-1}(q)\Si_n\bigl(1,(q-1)\Ss Q\bigr)\\&=\Si_d(q)\sum_{n<d+1}\Si_n(1)\Si_{n-1}^{-1}(q)\Si_n\bigl(1,(q-1)\bigr)\Si_{<n}(Q)\\&=-D_1^{-1}\Si_d(q)\sum_{n<d+1}\Si_n(1)\Si_{<n}(Q)\\&=\Si_{d+1}(1,q-1)\Si_{<d+1}(1,Q)\\&=\Si_{d+1}\bigl(1,(q-1)\Ss(1,Q)\bigr).
    \end{align*}
    To prove \Cref{eq:D1}, we have
    \begin{align*}
        \Si_d(q)\Si_{<d+1}(\prescript{1+}{}{P})&=\Si_d(q)\sum_{n<d+1}\Si_n(\prescript{1+}{}{P})\\&=\Si_d(q)\sum_{n<d+1}\Si_n(1)\Si_{n}(P)\\&=\Si_d(q)\sum_{n<d+1}\Si_n(1)\bigl(\Si_{<n+1}(P)-\Si_{<n}(P)\bigr)\\&=\Si_d(q)\sum_{n<d+1}\Si_n(1)\Si_n(q)^{-1}\Si_{n+1}\Bigl(\bigl(1,(q-1)\Ss Q\bigr)\Bigr)-\Si_{d+1}\bigl(1,(q-1)\Ss(1,Q)\bigr)\\&=-D_1^{-1}\Si_d(q)\sum_{n<d+1}\Si_n(1)\Si_{<n+1}(Q)-\Si_{d+1}\bigl(1,(q-1)\Ss(1,Q)\bigr)\\&=\Si_{d+1}(1,q-1)\sum_{n<d+1}\Si_n(1)\Si_{<n+1}(Q)-\Si_{d+1}\bigl(1,(q-1)\Ss(1,Q)\bigr)\\&=\Si_{d+1}(1,q-1)\bigl(\Si_{<d+1}(\prescript{1+}{}{Q})+\Si_{<d+1}(1,Q)\bigr)-\Si_{d+1}\bigl(1,(q-1)\Ss(1,Q)\bigr)\\&=\Si_{d+1}\bigl(1,(q-1)\Ss\prescript{1+}{}{Q}\bigr).
    \end{align*}
\end{proof}

\begin{theorem}\label{thm:P}
    The image of $\cP$ lies in $\br_{w+2q}$.
\end{theorem}
\begin{proof}
    By \cite{IKN2026} Proposition 6.1, we have binary relation 
    \begin{equation}\label{eq:basic}
        \Bigl[(q)(\St+\Sc)\bigl((2,q-1)+(1,1,q-1)\bigr);-\bigl(1,(q-1)\Ss(1,q)\bigr)\Bigl],
    \end{equation}
    which is actually $\cP\bigl((1)\bigr)$. In \Cref{lem:D}, set $P=(2,q-1)+(1,1,q-1)$ and $Q=(1,q)$, then $R=\cP\bigl((1)\bigr)$. Repeatedly apply $\cD_0$ and $\cD_1$ operators, note that $\cD_0\bigl(\cP(\fs)\bigr)=\cP\bigl((1,\fs)\bigr)$ and $\cD_1\bigl(\cP(\fs)\bigr)=\cP(\prescript{1+}{}{\fs})$, which implies our desired result.
\end{proof}

\subsection{Decomposition of $\cS_w^{(\Fq)}$ and $\ind$ operator}\label{subsec:dec}
For any $w\geq q$, and $\fs=(s_1,\cdots,s_l)\in\cI_w\setminus\BTTw$, write $\fs=(\fh,q+r,\ft)$, where $\fh=(s_1,\cdots,s_{k-1})\in[q]^{\bullet}$; either $r>0$, or $r=0$ and $\ft=\emptyset$ (in this case we set $(r,\ft)=\emptyset$). Set $\cL\bigl[\fh;(r,\ft)\bigr]=(\fh,q+r,\ft)+\mathbbm{1}_{r \neq 0}(\fh,q,r,\ft)$ and $\cH\bigl[\fh;(r,\ft)\bigr]=\mathbbm{1}_{\fh \neq \emptyset}\bigl(\fh^{+1},(q-1)\Ss(r,\ft)\bigr)+\bigl(\fh,1,(q-1)\Ss(r,\ft)\bigr)$. We may directly write $\cL[\fh;\fs]$ to denote $\cL\bigl[\fh;(s_1,\fs^T)\bigr]$, so does $\cH$. We view $\cL,\cH$ as the linear maps from $\Span_{\Fq}\{[\fs;\ft]\}$ to $\cS^{(\Fq)}$ by linear extension.

Im, Kim and Ngo Dac constructed the following binary relation (See \cite{IKN2026}, Section 4.2):
\[
    \begin{cases}
        \bigl[\cL[\fh;\ft];D_1\cH[\fh;\ft]\bigr], \quad&\text{if }\fh=\emptyset;\\
        \bigl[\cL[\fh;\ft]+D_1\cH[\fh;\ft];0\bigr], \quad&\text{if }\fh\neq\emptyset.
    \end{cases}
\]

To give a direct sum decomposition of $\cS_w^{(\Fq)}$, we first introduce the (depth, lex)-order which will be repeatedly used.

\begin{definition}
The (depth, lex)-order $\prec$ on $\cI$ is defined as follows.
For any two tuples $\fs, \ft \in \cI_w$, define $\fs \prec \ft$ if and only if
\begin{itemize}
\item either $\dep(\fs) < \dep(\ft)$,
\item or $\dep(\fs) = \dep(\ft)$ and $\fs$ is lexicographically smaller than $\ft$.
\end{itemize}
\end{definition}

\begin{lemma}\label{lem:decomposition}
    Set $\sL_w=\Span_{\Fq}\{\cL[\fh;\fs] ~|~ \fh\in[q]^{\bullet},~ \wt(\fh)+\wt(\fs)+q=w\}$. Then \[\cS_w^{(\Fq)}=\sL_w\oplus\Span_{\Fq}\BTTw.\] Moreover, $\cR_w$ is a subspace of $\sL_w$.
\end{lemma}
\begin{proof}
    For the direct sum decomposition, note that for any $\fs=(\fh,q+r,\ft)\in\cI_w\setminus\BTTw$, $\fs$ is exactly the smallest term among the tuples of $\cL\bigl[\fh;(r,\ft)\bigr]$ in the (depth, lex)-order. By comparing the smallest term, we obtain the direct composition.

    To prove $\cR_w\subseteq\sL_w$, for any $R\in\cR_w$, write $R=L+T$ with $L=\sum a_i\cL[\fh_i;\fs_i]\in\sL_w$ and $T\in\Span_{\Fq}\BTTw$. Let $H=\sum a_i\cH[\fh_i;\fs_i]$. Since $R\in\cR_w$, we have $\Li(R)=\Li(L)+\Li(T)=0$. But $\Li(L)=-D_1\Li(H)$, we have $\Li(T)-D_1\Li(H)=0$. Now repeat this procedure for $H$ until all tuples belonging to $\BTTw$, then we obtain a $\Fq(D_1)$-linear relation among $\BTTw$ with the zero-degree part $T$. By the \Cref{thm_CCMIKLNP}, there have to be $T=0$.
\end{proof}
\begin{remark}
    The proof of $\cR_w\subseteq\sL_w$ is essentially the one in \cite{IKN2026}.
\end{remark}

Next we define $\ind:\sL_w\to\cS_w^{(\Fq)}$ by $\ind\bigl(\cL[\fh;\fs]\bigr)=\cH[\fh;\fs]$ for $\fh\in[q]^{\bullet}$, which is exactly the $\ind$ operator in \cite{IKN2026} but is extended to $\sL_w$. It is easy to see $\ind$ preserves $\cR_w$ by definition. Moreover, it is nilpotent on $\cR_w$.

\begin{lemma}\label{lem:nil}
    For any $w\geq q$, $\ind$ acts on $\cR_w$ nilpotently.
\end{lemma}
\begin{proof}
    For any $R\in\cR_w$, note that if $\ind(R)\neq0$, then in the (depth, lex)-order $\prec$, $\min_{\prec}R\prec\min_{\prec}\ind(R)$, where $\min_{\prec}P$ means the smallest tuple of $P$. Since there are only finitely many tuples for a fixed weight, there exists a $k\in\bbZ+$ such that $\ind^k=0$.
\end{proof}

The following theorem will be used to prove our final result.
\begin{theorem}\label{thm:indinverse}
    For any $w\geq q$ and subspace $L\subseteq\cR_w$, if $\ind^{-1}L\subseteq L$, then $L=\cR_w$.
\end{theorem}
\begin{proof}
    It is a direct corollary of \Cref{lem:nil}.
\end{proof}

\subsection{Kernel of $\cL$}\label{subsec:KernelL}
We now determine $\ker\cL$, where we can see the similarity in spirit between the method in this paper and the idea in \cite{HL2026}.

For any $w\geq q$, $\fs=(s_1,\cdots,s_l)\in\cI_w\setminus\BTTw$, and $s_i$ which is a component of $\fs$ such that $s_i>q$ or $s_i=q$ and $i=l$ (Here we don't require $s_i$ is the first component satisfying this condition), write $\fs=(\fh_i,q+r_i,\ft_i)$, $\fs^i=(\fh_i,q,r_i,\ft_i)$ and $E_i(\fs)=\bigl[\fh_i;(r_i,\ft_i)\bigr]$. Use the previous notation, we still write $\fs=(\fh,q+r,\ft)$ if $q+r$ is the first component satisfying condition, and we call the position  of $q+r$ is a canonical position. We also set $E(\fs)=\bigl[\fh;(r,\ft)\bigr]$, $\fs^{cano}=(\fh,q,r,\ft)$.

\begin{theorem}
    For any $w\in\bbZ+$, $\sL_w=\Span_{\Fq}\{\cL[\fh;\fs] ~|~ \wt(\fh)+\wt(\fs)+q=w\}$
\end{theorem}
\begin{proof}
    Denote the right hand side by $\sL_w'$. For $\fs=(s_1,s_2,\cdots,s_l)\in\cI_w$, it can be uniquely written as $(qk_1+r_1,qk_2+r_2,\cdots,qk_l+r_l)$ with $k_i\geq0,1\leq r_i\leq q$. Define the $\Fq$-linear map $\cE:\cS_w^{(\Fq)}\to\Span_{\Fq}\BTTw$ by
    \[
    \cE(\fs)=\begin{cases}
        (-1)^{\sum k_i}(\underbrace{q,\cdots,q}_{k_1\text{ q's}},r_1,\underbrace{q,\cdots,q}_{k_2\text{ q's}},r_2,\cdots,\underbrace{q,\cdots,q}_{k_l\text{ q's}},r_l), &r_l\neq q,\\
         0, &r_l=q
    \end{cases}
    \]

    It is easy to see $\cE|_{\Span_{\Fq}\BTTw}=\operatorname{id}$ and $\cE\circ\cL=0$.

    Define another $\Fq$-linear map $f:\cS^{(\Fq)}\to\Span_{\Fq}\{[\fs;\ft]\}$ inductively in the reverse (dep,lex)-order \[f(\fs)=\begin{cases}
    0, &\fs\in\BTTw,\\
    E(\fs)-f(\fs^{cano}), &\fs\in\cI_w\setminus\BTTw,\fs^{cano}\neq\fs,\\
    E(\fs),&\fs\in\cI_w\setminus\BTTw,\fs^{cano}=\fs.
    \end{cases}\]

    Now we can begin the proof, by induction in the reverse (dep,lex)-order, we have $\cL\circ f=\operatorname{id}-\cE$, which implies $\ker\cE=\{P\in\cS_w^{(\Fq)} ~|~ P=\cL\circ f(P)\}$. Note that the image of $f$ lies in $\Span_{\Fq}\{[\fh;\fs] ~|~ \fh\in[q]^{\bullet},~ \wt(\fh)+\wt(\fs)+q=w\}$, we have $\ker \cE\subseteq\sL_w$. Since $\cE\circ\cL=0$, we have $\ker \cE=\sL_w$. Because $\sL_w'\subseteq\ker\cE=\sL_w\subseteq\sL_w'$, we are done.
\end{proof}

\begin{theorem}\label{thm:kerL}
    For any $w\in\bbZ+$, Set \[Q_w=\Span_{\Fq}\{-E_i(\fs)+E_j(\fs)+E_{j+1}(\fs^i)-\mathbbm{1}_{r_j\neq0}E_i(\fs^j) ~|~ \fs\in\cI_w,~i<j\text{ run over all possible positions}\},\]
    then $\ker\cL=Q_w$.
\end{theorem}
\begin{proof}
    We first prove that for any $\fs\in\cI_w$ and possible $i$, we have $E_i(\fs)\equiv f\circ\cL\bigl(E_i(\fs)\bigr)\mod{Q_w}$.
    
    We do it by induction on $\fs$ with respect to reverse (dep,lex)-order.
    
    When $E_i=E$ (i.e. $i$ is canonical), it is an equality. In other cases, if $r_i>0$, we have \begin{equation}\label{eq:congruent}
        -E_i(\fs)+E(\fs)\equiv-E(\fs^i)+E_{i+1}(\fs^{cano}).
    \end{equation}
    We first show $-E(\fs^i)+E_{i+1}(\fs^{cano})\equiv -f(\fs^i)+f(\fs^{cano})$. It suffices to show \[-E_{i+1}(\fs^{cano})+E(\fs^{cano})\equiv f\bigl((\fs^{cano})^{cano}\bigr)-f\bigl((\fs^i)^{cano}\bigr),\] Using \Cref{eq:congruent}, it suffices to show \[-E\bigl((\fs^{cano})^{(i+1)}\bigr)+E_{i+2}\bigl((\fs^{cano})^{cano}\bigr)\equiv f\bigl((\fs^{cano})^{cano}\bigr)-f\bigl((\fs^i)^{cano}\bigr).\]
    Note that $(\fs^i)^{cano}=(\fs^{cano})^{i+1}$, inductive assumption implies desired result.

    If $r_i=0$, then \[E_i(\fs)\equiv E(\fs)-E_{i+1}(\fs^{cano}).\]

    Inductive assumption implies
    \[E_{i+1}(\fs^{cano})\equiv f(\fs^{cano}),\]
    hence $E_i(\fs)\equiv f(\fs).$

    Then we have \[E_i(\fs)\equiv f(\fs)+\mathbbm{1}_{r_i>0}f(\fs^i)\equiv f\circ\cL\bigl(E_i(\fs)\bigr).\]

    Now for any $R\in\ker\cL$, we can always find some $\fs_i\in\cS_w^{(\Fq)}$ such that $R=\sum a_iE_{j_i}(\fs_i)$, hence $R\equiv 0$. So $\ker\cL\subseteq Q_w$, the converse can be checked directly.
\end{proof}

\subsection{General $\cP$ operator and cut identities}
In ideas and techniques, this subsection is very similar to \cite{Hu2026}.
The goal of this subsection is to construct a more general $\cP$ operator, and realize it by using the original $\cC,\cP$ operators. For $\fs,\ft\in\cI$, we define
\[\cP[\fs;\ft]=\biggl[(q)(\St+\Sc)\bigl(\fs^{+1}+(\fs,1),(q-1)\Ss\ft\bigr);-\Bigl(1,(q-1)\Ss\bigl(\fs,(q)(\St+\Sc)\ft\bigr)\Bigr)\biggr].\]
Note that it is the previous $\cP$ operator if set $\ft=\emptyset$. We will show the image of general $\cP$ still lies in $\br$.

To do this, we first introduce the cut of tuples. For a $\fs=(s_1,\cdots,s_l)\in\cI$, a cut of $\fs$ is a pair of tuple $\fu,\fv$, such that $\fs=(\fu,\fv)$, i.e. $\fu=(s_1,\cdots,s_k),\fv=(s_{k+1},\cdots,s_l)$ for some $k$, possibly $\fu=\emptyset$ or $\fv=\emptyset$. The following cut identities are complex but intuitive.

\begin{lemma} For any $\fs,\ft\in\cI\setminus\{\emptyset\}$ and $c\in\bbZ+$,
    \begin{enumerate}
        \item \begin{equation}\label{eq:cut1}
            (\fs,c)\Ss\ft=\sum_{\ft=(\fu,\fv)}\bigl(\fs\Ss\fu,(c)(\St+\Sc)\fv\bigr);
        \end{equation}
        \item \[\label{eq:cut2}
            \bigl(\fs^{+1}+(\fs,1)\bigr)\Ss\ft=\sum\limits_{\ft=(\fu,\fv)}\bigl((\fs\Ss\fu)^{+1}+(\fs\Ss\fu,1),\fv\bigr);
        \]
        \item \begin{equation}\label{eq:cut3}
            \bigl(\fs^{+1}+(\fs,1),q-1\bigr)\Ss\ft=\sum\limits_{\ft=(\fu,\fv)}\bigl((\fs\Ss\fu)^{+1}+(\fs\Ss\fu,1),(q-1)\Ss\fv\bigr);
        \end{equation}
    \end{enumerate}
    where we set $(c)(\St+\Sc)\emptyset=(c)$
\end{lemma}
\begin{proof}
    For $\fa=(a_1,\cdots,a_k),\fb=(b_1,\cdots,b_l)$, \Cref{eq:Ss} yields that
    \begin{equation}\label{eq:extendss}
        \fa\Ss\fb=\sideset{}{'}{\sum}(c_1|c_2|\cdots|c_{k+l})
    \end{equation}
    where the restricted sum means: \begin{enumerate}
        \item[(1)] $c_1c_2\cdots c_{k+l}$ is a permutation of $a_1,\cdots,a_k,b_1,\cdots,b_l$ preserving the order of $a_i$'s and $b_i$'s;
        \item[(2)] ``|'' is ``+'' or ``,'' satisfying there are at most one $a_i$ and one $b_i$ between a pair of adjacent ``,'';
        \item[(3)] if two terms differ only by a commutation of addition, they are counted once in the sum (for example, $(a_1+b_1,a_2,b_2)$ and $(b_1+a_1,a_2,b_2)$ are only counted once).
    \end{enumerate} 

    Then (i) is directly yielded by \Cref{eq:extendss}.

    For (ii), use \Cref{eq:extendss} to expand the left hand side, there are terms of forms $(c_1|c_2|\cdots|a_k+1|\cdots|c_{k+l})$ corresponding to $\fs^{+1}\Ss\ft$ and $(c_1|c_2|\cdots|a_k|\cdots|1|\cdots|c_{k+l})$ corresponding to $(\fs,1)\Ss\ft$. By (3), we may assume that the $(c_1|c_2|\cdots|a_k+1|\cdots|c_{k+l})=(c_1|c_2|\cdots|a_k+1,\cdots|c_{k+l})$ and $(c_1|c_2|\cdots|a_k|\cdots|1|\cdots|c_{k+l})=(c_1|c_2|\cdots|a_k|\cdots|1,\cdots|c_{k+l})$. Set $\fu=(b_1,\cdots,b_n)$ such that $b_n$ precedes $a_l+1$ or $1$, and $b_{n+1}$ follows $a_l+1$ or $1$,. Summing over the terms with the same $\fu$, it is $\bigl((\fs\Ss\fu)^{+1}+(\fs\Ss\fu,1),\fv\bigr)$, so we are done.

    A similar argument to (ii) proves (iii).
\end{proof}

\begin{theorem}\label{thm:CP}
    For any $\fs,\ft\in\cI$, $\cC_{\ft}\cP(\fs)=\sum\limits_{\ft=(\fu,\fv)}\cP[\fs\Ss\fu;\fv]$
\end{theorem}
\begin{proof}
    It suffices to show \[\biggl[\Bigl((q)(\St+\Sc)\bigl(\fs^{+1}+(\fs,1),q-1\bigr)\Bigr)(\St+\Sc)\ft;-\bigl(1,(q-1)\Ss(\fs,q)\bigr)\St\ft\biggr]=\sum\limits_{\ft=(\fu,\fv)}\cP[\fs\Ss\fu;\fv]\]

    It is easy to check that for any $P,Q,S\in\cS^{(\Fq)}$, $\bigl(P(\St+\Sc)Q\bigr)(\St+\Sc)S=P(\St+\Sc)(Q\Ss S)$ (check it hold for single tuples then extend linearly). So we need to show
    \[\biggl[(q)(\St+\Sc)\Bigl(\bigl(\fs^{+1}+(\fs,1),q-1\bigr)\Ss\ft\Bigr);-\bigl(1,(q-1)\Ss(\fs,q)\bigr)\St\ft\biggr]=\sum\limits_{\ft=(\fu,\fv)}\cP[\fs\Ss\fu;\fv]\]\

    Using \Cref{eq:cut3}, we have
    \[\text{LHS}=\biggl[\sum\limits_{\ft=(\fu,\fv)}(q)(\St+\Sc)\Bigl(\bigl((\fs\Ss\fu)^{+1}+(\fs\Ss\fu,1),(q-1)\Ss\fv\bigr)\Bigr);-\bigl(1,(q-1)\Ss(\fs,q)\bigr)\St\ft\biggr]
    \]

    Using \Cref{eq:cut1},
    \begin{align*}
        \text{RHS}&=\sum\limits_{\ft=(\fu,\fv)}\biggl[(q)(\St+\Sc)\bigl((\fs\Ss\fu)^{+1}+(\fs\Ss\fu,1),(q-1)\Ss\fv\bigr);-\Bigl(1,(q-1)\Ss\bigl(\fs\Ss\fu,(q)(\St+\Sc)\fv\bigr)\Bigr)\biggr]\\&=\biggl[\sum\limits_{\ft=(\fu,\fv)}(q)(\St+\Sc)\Bigl(\bigl((\fs\Ss\fu)^{+1}+(\fs\Ss\fu,1),(q-1)\Ss\fv\bigr)\Bigr);-\bigl(1,(q-1)\Ss(\fs,q)\Ss\ft\bigr)\biggr]
    \end{align*}

    Note that $\bigl(1,(q-1)\Ss(\fs,q)\bigr)\St\ft=\bigl(1,(q-1)\Ss(\fs,q)\Ss\ft\bigr)$, we done.
\end{proof}

Now we show the general $\cP$ can be realized by using the original $\cC$ and $\cP$ operators.
\begin{corollary}\label{cor:P}
    For any $\fs,\ft\in\cI$, $\cP[\fs;\ft]$ is a linear combination of $\cC_{\fa}\cP(\fb)$ for some $\fa,\fb$.
\end{corollary}
\begin{proof}
    We prove it by induction on $\dep(\ft)$. When $\dep(\ft)=0$, $\ft=\emptyset$, there is nothing to prove. Assume we have done for all $\ft$ with $\dep<d$. Now for a $\ft$ with depth $d$, by previous theorem, we have
    \[\cP[\fs;\ft]=\cC_{\ft}\cP(\fs)-\sum_{\ft=(\fu,\fv),\fu\neq\emptyset}\cP[\fs\Ss\fu;\fv].\]

    The $\fv$'s in RHS all have depth less than $d$. By inductive assumption, RHS can be written as a linear combination of $\cC_{\fa}\cP(\fb)$, then we are done.
\end{proof}

We now use \Cref{cor:P} to describe $\BCP_w$.

\begin{theorem}\label{thm:BCPHL}
    For any $w\in\bbZ_+$, $\BCP_w=\cH(\ker\cL)$.
\end{theorem}
\begin{proof}
    We first show $\cH(\ker\cL)\subseteq\BCP_w$. By \Cref{thm:kerL}, it suffices to show $\cH\bigl(-E_i(\fs)+E_j(\fs)+E_{j+1}(\fs^i)-\mathbbm{1}_{r_j\neq0}E_i(\fs^j)\bigr)\in\BCP_w$ for any $\fs\in\cI_w$ and possibly $i,j$.

    For a $\fs=(\fh,q+r_i,\fm,q+r_j,\ft)$ with $r_i,r_j>0$, we have
    \[
        -E_i(\fs)+E_j(\fs)+E_{j+1}(\fs^i)-E_i(\fs^j)=-\Bigl[\fh;\bigl(r_i,\fm,(q)(\St+\Sc)(r_j,\ft)\bigr)\Bigr]+\Bigl[\bigl(\fh,(q)(\St+\Sc)(r_i,\fm)\bigr);(r_j,\ft)\Bigr].
    \]
    Then
    \begin{align*}
        \cH\bigl(-E_i(\fs)+E_j(\fs)+E_{j+1}(\fs^i)-E_i(\fs^j)\bigr)=&-\Bigl(\fh^{+1}+(\fh,1),(q-1)\Ss\bigl(r_i,\fm,(q)(\St+\Sc)(r_j,\ft)\bigr)\Bigr)\\&+\biggl(\fh,\Bigl((q)(\St+\Sc)\bigl((r_i,\fm)^{+1}+(r_i,\fm,1)\bigr),(q-1)\Ss(r_j,\ft)\Bigr)\biggr)
    \end{align*}

    We can see
    \[\cP\bigl[(r_i,\fm);(r_j,\ft)\bigr]=\biggl[(q)(\St+\Sc)\bigl((r_i,\fm)^{+1}+(r_i,\fm,1),(q-1)\Ss(r_j,\ft)\bigr);-\Bigl(1,(q-1)\Ss\bigl(r_i,\fm,(q)(\St+\Sc)(r_j,\ft)\bigr)\Bigr)\biggl].\]

    So we have $\cBs_{\fh}\Bigl(\cP\bigl[(r_i,\fm);(r_j,\ft)\bigr]\Bigr)=\Bigl[\cH\bigl(-E_i(\fs)+E_j(\fs)+E_{j+1}(\fs^i)-E_i(\fs^j)\bigr);0\Bigr]$ (if $\fh=\emptyset$, we need a lift).

    For $\fs=(\fh,q+r_i,\fm,q)$ with $r_i>0$, we can obtain desired results similarly.

    For $\BCP_w\subseteq\cH(\ker\cL)$, it suffices to show $\overline{\cBs_{\fh}\bigl(\cP[\fs;\ft]\bigr)}\in\cH(\ker\cL)$ for all $\fh,\fs,\ft\in\cI$. Let $R=\Bigl[\bigl(\fh,\cL[\emptyset;\fs]\bigr);\ft\Bigr]-\bigl[\fh;\cL[\fs;\ft]\bigr]$, direct calculation shows $\cL(R)=0$ and $\cH(R)=\overline{\cBs_{\fh}\bigl(\cP[\fs;\ft]\bigr)}$, so we are done.
\end{proof}

\begin{remark}
    If we admit the results in \cite{HL2026}, then we have done, because $\cH(\ker\cL)$ is actually the space spanned by all quadruple-carry relations.
\end{remark}

\subsection{Separate singletons}
Let $\sC=\Span_{\Fq}\{[\fh;\ft] ~|~ \fh\neq\emptyset\}$. We study $\im\cH|_{\sC}$ precisely in this subsection. 

It is easy to see that
\begin{equation}\label{eq:H+}
    \cH(\sC)=\Span_{\Fq}\{\bigl(\fs,(a+1)+(a,1),(q-1)\Ss\ft\bigr) ~|~ \fs,\ft\in\cI,a\in\bbZ+\}.
\end{equation}

Denote the above set by $H_+$, we will show the following lemma.

\begin{lemma}\label{lem:separate}
    \begin{enumerate} For any $P\in\cS^{(\Fq)}$
        \item $P\in H_+\iff(1,P)\in H_+$;
        \item $P\in H_+\iff(q-1)\Ss P\in H_+$
    \end{enumerate}
\end{lemma}
\begin{proof}
    \begin{enumerate}
        \item ``$\implies$'' is obvious by \Cref{eq:H+}, in fact $H_+$ is stable under adding any prefix (we will repeatedly use this fact).

        To prove ``$\impliedby$'', for a $A\in H_+$, using \Cref{eq:H+}, it can be written as
        \[
            A=\sum_{i}\bigl((i+1)+(i,1),B_i\bigr)+\sum_{j}(j,C_j),\quad\text{for some } i,j\in\bbZ_+, B_i\in(q-1)\Ss\cS^{(\Fq)},C_j\in H_+.
        \]
        To see this, just write $A$ as a linear combination of elements in RHS of \Cref{eq:H+}, then separate the terms with $\fs=\emptyset$ and the terms with $\fs\neq\emptyset$, for the latter ones, replace $\fs$ by $\fs^T$ to obtain $C_j$.

        Now we define the extraction operator $\rho_n$ for $n\in\bbZ_+$. For a tuple $\fs=(s_1,\fs^T)$, we define
        \[
            \rho_n(\fs)=\begin{cases}
                0, &s_1\neq n\\
                \fs^T, &s_1=n
            \end{cases}
        \]

        Apply it to $A$, we have
        \[
            \rho_1(A)=(1,B_1)+C_1,
        \]
        \[
            \rho_n(A)=B_{n-1}+(1,B_n)+C_n, \quad n\geq2.
        \]
        
        There is a smallest $N$ such that $\rho_n(A)=0$ for any $n>N$. So for $n\geq N$, we have $B_n=-(1,B_{n+1})-C_{n+1}$. For $n$ large enough, $B_n=0$, hence $B_n\in H_+$. Since $C_n\in H_+$ and $B_{n+1}\in H_+\implies(1,B_{n+1})\in H_+$, we have $B_n\in H_+$ for all $n\geq N$ by backward induction. As result, we have
        \begin{equation}\label{eq:N1}
            \rho_1(A)=(1,B_1)+C_1\in H_+\quad\text{if }N=1,
        \end{equation}
        \begin{equation}\label{eq:N2}
            \rho_N(A)=B_{N-1}+(1,B_N)+C_N\in (q-1)\Ss\cS^{(\Fq)}+H_+\quad\text{if }N\geq2.
        \end{equation}

        Apply \Cref{eq:N1} to $A=(1,P)$, we have $P\in H_+$.

        \item To prove ``$\implies$'', note that
        \[(q-1)\Ss(\fs,\ft)=\bigl((q-1)\Ss\fs,\ft\bigr)+\bigl(\fs,(q-1)\Ss\ft\bigr)-(\fs,q-1,\ft).\]
        So we have
        \begin{align*}
            (q-1)\Ss\bigl(\fs,(a+1)+(a,1),(q-1)\Ss\ft\bigr)=&\bigl((q-1)\Ss\fs,(a+1)+(a,1),(q-1)\Ss\ft\bigr)\\&+\Bigl(\fs,(q-1)\Ss\bigl((a+1)+(a,1),(q-1)\Ss\ft\bigr)\Bigr)\\&-\bigl(\fs,q-1,(a+1)+(a,1),(q-1)\Ss\ft\bigr),
        \end{align*}

        where the first and the last terms in RHS belong to $H_+$, so we only need to prove \[\Bigl(\fs,(q-1)\Ss\bigl((a+1)+(a,1),(q-1)\Ss\ft\bigr)\Bigr)\in H_+.\]
        It can be implied from the following calculation
        \begin{align*}
            (q-1)&\Ss\bigl((a+1)+(a,1),(q-1)\Ss\ft\bigr)\\=&\bigl((q+a)+(q+a-1,1),(q-1)\Ss\ft\bigr)+\bigl(q-1,(a+1)+(a,1),(q-1)\Ss\ft\bigr)\\&+\bigl((a,q)+(a,q-1,1),(q-1)\Ss\ft\bigr)+\bigl((a+1)+(a,1),(q-1)\Ss(q-1)\Ss\ft\bigr).
        \end{align*}

        To show ``$\impliedby$'', we first show that if $(q-1)\Ss(1,P)\in H_+$, then $(q-1)\Ss P\in H_+$. In fact, the greatest first component of the tuples in $(q-1)\Ss(1,P)$ is $q$, so $\rho_q\bigl((q-1)\Ss(1,P)\bigr)\in(q-1)\Ss\cS^{(\Fq)}+H_+$ by \Cref{eq:N2}. Since $\rho_q\bigl((q-1)\Ss(1,P)\bigr)=P$, we have $\bigl((q)+(q-1,1),P\bigr)\in H_+$. Since $(q-1)\Ss(1,P)=(q,P)+(q-1,1,P)+\bigl(1,(q-1)\Ss P\bigr)$, we are done by using (i).

        Now we prove $(q-1)\Ss P\in H_+\implies P\in H_+$ by induction on $\wt(P)$. For $\wt(P)=0$, there is nothing to prove. Assume we have done for $\wt<w$, now for $P\in\cS_w^{(\Fq)}$, write it as $P=\sum\limits_{1\leq i\leq m}(i,P_i)$, where $m$ is the maximal first component among the tuples of $P$, then \[(q-1)\Ss P=(q-1,P)+\sum_{1\leq i\leq m
        }(q-1+i,P_i)+\sum_{1\leq i\leq m
        }\bigl(i,(q-1)\Ss P_i\bigr)\]
        Use \Cref{eq:N2}, we have \[\rho_{q-1+m}\bigl((q-1)\Ss P\bigr)=P_m\in(q-1)\Ss\cS^{(\Fq)}+H_+.\]
        We may choose a decomposition $P_m=S+T$ with $S\in(q-1)\Ss\cS^{(\Fq)}, T\in H_+$, then \[P=\sum_{1\leq i\leq m-1}(i,P_i)+(m,S)+(m,T).\]
        Since $(m-1,1,S)+(m,S)\in H_+$, we have
        \[P\equiv\sum_{1\leq i\leq m-1}(i,P_i)-(m-1,1,S)\mod{H_+}.\]
        Repeatedly use the above argument, we have
        \[P\equiv(1,Q)\mod{H_+},\quad\text{for some }Q\in \cS^{(\Fq)}.\]
        Since $(q-1)\Ss H_+\subseteq H_+$, we have $(q-1)\Ss(1,Q)\in H_+$, which implies $(q-1)\Ss Q\in H_+$. As $\wt(Q)<\wt(P)$, we have $Q\in H_+$. Then by (i), we are done.
    \end{enumerate}
\end{proof}

\subsection{Kernel of $\cH$}
The final of our preparation is determining $\ker\cH$.

Define $\Fq$-linear map $\Lambda:\sC\to\sC$ by
\[\Lambda[\fs;\ft]=\bigl[(q-1)\Ss\fs;\ft\bigr]+\bigl[\fs;(q-1)\Ss\ft\bigr].\]

\begin{lemma}\label{lem:Lambda}
    For any $[P;Q]\in\sC$,
    \begin{enumerate}
        \item $\cL\bigl(\Lambda[P;Q]\bigr)=(q-1)\Ss\cL[P;Q]$;
        \item $\cH\bigl(\Lambda[P;Q]\bigr)=(q-1)\Ss\cH[P;Q]$.
    \end{enumerate}
\end{lemma}
\begin{proof}
It suffices to show theorem for $[P;Q]=[\fs;\ft]$.
    \begin{enumerate}
        \item We have \[\text{LHS}=\bigl((q-1)\Ss \fs,\prescript{q+}{}{\ft}\bigr)+\bigl((q-1)\Ss\fs,q,\ft\bigr)+\Bigl(\fs,\prescript{q+}{}{\bigl((q-1)\Ss\ft\bigr)}\Bigr)+\bigl(\fs,q,(q-1)\Ss\ft\bigr),\]
        \[\text{RHS}=(q-1)\Ss\bigl((\fs,\prescript{q+}{}{\ft})+(\fs,q,\ft)\bigr).\]
        Note that \[\prescript{q+}{}{\bigl((q-1)\Ss\ft\bigr)}=(q-1)\Ss\prescript{q+}{}{\ft}-(q-1,\prescript{q+}{}{\ft})+(2q-1,\ft),\]
        then directly calculation implies.

        \item Note that \[\bigl((q-1)\Ss\fs\bigr)^{+1}=(q-1)\Ss\fs^{+1}+(\fs,q)-(\fs^{+1},q-1),\]
        Similar calculation implies.
    \end{enumerate}
\end{proof}

\begin{theorem}\label{thm:KerH}
    \[\ker\cH=\Span_{\Fq}\{\bigl[\fs;\cH(R)\bigr]-\bigl(\fs^{+1}+(\fs,1),\Lambda(R)\bigr)~|~R\in\sC\},\]
    where $\bigl(\fs,[\fu;\fv]\bigr)=\bigl[(\fs,\fu);\fv\bigr]$, and we set $\emptyset^{+1}+(\emptyset,1)=(1)$.
\end{theorem}
\begin{proof}
    Using,
    \[\cH(\fh,R)=\bigl(\fh,\cH(R)\bigr),\quad\text{for any }\fh\in \cI,R\in\sC.\]
    Direct calculation implies ``$\supseteq$''.
    
    We prove ``$\subseteq$'' by induction on $\wt(P)+\wt(Q)$ for $\cH[P;Q]=0$. If $\wt(P)+\wt(Q)=0$, then $[P;Q]=k[\emptyset;\emptyset]$ for some $k\in\Fq$, the result is obvious. Assume we have done for $\wt(P)+\wt(Q)<w$. For a $R\in\ker\cH$, it can be written as a sum $\sum[\emptyset;T]+\sum[P;Q]$ for some $T\in\cS^{(\Fq)}$ and $[P;Q]\in\sC$, then
    \[0=\cH(R)=\sum\bigl(1,(q-1)\Ss T\bigr)+\sum\bigl(P^{+1}+(P,1),(q-1)\Ss Q\bigr).\]
    The second term in RHS lies in $H_+$, hence $\bigl(1,(q-1)\Ss T\bigr)\in H_+$.
    
    \Cref{lem:separate} implies $T\in H_+$, so we can choose $R_1\in\sC$ such that $\cH(R_1)=T$. Replace $R$ by $R-\Bigl(\bigl[\emptyset;\cH(R_1)\bigr]-\bigl(1,\Lambda(R_1)\bigr)\Bigr)$, we may assume $R\in\ker\cH\cap\sC$. For such $R\in\ker\cH\cap\sC$, it can be written as \[R=\sum\bigl[(n);P_n\bigr]+\sum(m,S_m),\quad P_n\in\cS^{(\Fq)},S_m\in\sC,\]
    where we separate terms depending the depths of the first part.
    
    Then apply $\cH$ to $R$, we have
    \[0=\sum\bigl((n+1)+(n,1),(q-1)\Ss P_n\bigr)+\sum\bigl(m,\cH(S_m)\bigr).\]
    Use the similar argument in the proof of \Cref{lem:separate} (ii), we can show $(q-1)\Ss P_n\in H_+$ for all $n$. \Cref{lem:separate} implies $P_n\in H_+$, so there are some $R_n\in\sC$ such that $\cH(R_n)=P_n$.

    Replace $R$ by $R-\sum\Bigl(\bigl[(n);\cH(R_n)\bigr]-\bigl((n+1)+(n,1),\Lambda(R_n)\bigr)\Bigr)$, then $R$ has form of $\sum(m,S_m)$ with $S_m\in\sC$. Since the depths of their first part are at least $2$, the first components of $\cH(m,S_m)$ have to be different, which implies $\cH(m,S_m)=0$ for all $m$. Note that $\cH(m,S_m)=\bigl(m,\cH(S_m)\bigr)$, by the inductive assumption, we are done.
\end{proof}

\subsection{Proof of \Cref{thm:BCP}}
We can prove \Cref{thm:BCP} now. We first prove the final lemma in this section.
\begin{lemma}\label{lem:LBP} Let $R=[\fs;\ft]\in\sC$
    \[\cL\Bigl(\bigl[\fu;\cH(R)\bigr]-\bigl(\fu^{+1}+(\fu,1),\Lambda(R)\bigr)\Bigr)=\overline{\cBs_\fu\bigl(\cP(R)\bigr)},\]
    where we set $\cBs_{\emptyset}=\operatorname{id}$ for convenience.
\end{lemma}
\begin{proof}
    A complicated but straightforward calculations shows
    \begin{align*}
        \text{RHS}=\Bigl(\fu,\prescript{q+}{}{\bigl(\fs^{+1}+(\fs,1),(q-1)\Ss\ft\bigr)+\bigl(q,\fs^{+1}+(\fs,1),(q-1)\Ss\ft\bigr)}\Bigr)&\\-\Bigl(\fu^{+1}+(\fu,1),(q-1)\Ss\bigl(\fs,\prescript{q+}{}{\ft}+(q,\ft)\bigr)\Bigr)&.
    \end{align*}

    LHS needs a little attention:
    \[\text{LHS}=\Bigl(\fu,\prescript{q+}{}{\cH(R)}+\bigl(q,\cH(R)\bigr)\Bigr)-\bigl(\fu^{+1}+(\fu,1),(q-1)\Ss\cL(R)\bigr),\]
    where the right first term is direct, the second term needs \Cref{lem:Lambda} and the fact
    \[\cL(\fh,R)=\bigl(\fh,\cL(R)\bigr),\quad\text{for any }\fh\in \cI,R\in\sC.\]

    By the definition of $\cH,\cL$, we are done.
\end{proof}

Combine \Cref{cor:P}, \Cref{thm:KerH} and \Cref{lem:LBP}, we have $\cL(\ker\cH)\subseteq\BCP$.

Given an $P\in\sL_w$, there exists a $R=\sum a_i[\fh_i;\fs_i]$ with $\fh_i\in[q]^{\bullet}$, $\wt(\fh_i)+\wt(\fs_i)+q=w$ such that $P=\cL(R)$. Then $\cH(R)=\ind P$. 
If $\ind P=\cH(R)\in\BCP$, by \Cref{thm:BCPHL}, we have $R\in\ker\cL+\ker\cH$, which implies $P\in\cL(\ker\cH)$. So we get $\ind^{-1}(\BCP)\subseteq\cL(\ker\cH)$. So we have \[\ind^{-1}(\BCP)\subseteq\cL(\ker\cH)\subseteq\BCP.\]
Now use \Cref{thm:indinverse}, we prove \Cref{thm:BCP}.

\section{Dimension of $\cR_w$}\label{sec:dimension}
We determine the dimension of $\cR_w$, hence $\cZ_w$. Our method is relatively simpler than that in \cite{HL2026}.

\subsection{$\Fq+\Fq D_1$-relation}
Our starting point is the following set
\[U_w=\{P+D_1Q~|~P,Q\in\cS_w^{(\Fq)}~,~\Li(P+D_1Q)=0\},\]

\begin{lemma}\label{lem:U} For any $w\in\bbZ_+$
    \[U_w=\Span_{\Fq}\{\overline{\cBs_{\fh}\cC_{\ft}(R_1)}~|~\wt(\fh)+\wt(\ft)+q=w\}\]
\end{lemma}
\begin{proof}
    Direct calculation shows
    \[\overline{\cBs_{\fh}\cC_{\ft}(R_1)}=\bigl(\fh,\prescript{q+}{}{\ft}+(q,\ft)\bigr)+D_1\bigl(\fh^{+1}+(\fh,1),(q-1)\Ss\ft\bigr)\in U_w,\]
    where we set $\prescript{q+}{}{\emptyset}+(q,\emptyset)=(q)$ for convenience.

    For reverse inclusion, write $P=\sum\limits_{\fu_i\notin\BTTw}a_i\fu_i+\sum\limits_{\fv_j\in\BTTw}b_j\fv_j$. Note that each $\fu_i\notin\BTTw$, it can be written as $(\fh_i,\prescript{q+}{}{\ft_i})$ for some $\fh_i,\ft_i$. Repeatedly subtracting these $a_i\overline{\cBs_{\fh_i}\cC_{\ft_i}(R_1)}$, we can eliminate all terms not in $\BTTw$, with \[\Li\bigl(P+D_1Q-\sum a_i\overline{\cBs_{\fh_i}\cC_{\ft_i}(R_1)}\bigr)=0.\]
    Set $P+D_1Q-\sum a_i\overline{\cBs_{\fh_i}\cC_{\ft_i}(R_1)}=P'+D_1Q'$, then $P'\in\Span_{\Fq}\BTTw$. By comparing the degree of $D_1$ and using \Cref{thm_CCMIKLNP}, we have $P'=0$. So it suffices to show that if $Q\in\cR_w$, then $D_1Q\in\Span_{\Fq}\{\overline{\cBs_{\fh}\cC_{\ft}(R_1)}~|~\wt(\fh)+\wt(\ft)+q=w\}$.

    Direct calculation shows that 
    \[D_1\cP[\fs;\ft]=\cBs_{\cL[\emptyset;\fs]}\cC_{\ft}(R_1)-\cC_{\bigl(\fs,\cL[\emptyset;\ft]\bigr)}(R_1),\]
    then \Cref{thm:BCP} and \Cref{cor:P} imply the desired result.
\end{proof}

\begin{corollary}\label{cor:dimU}
    For any $w\in\bbZ_+$, $\dim_{\Fq}U_w=\dim_{\Fq}\sL_w+\dim_{\Fq}\cR_w.$
\end{corollary}
\begin{proof}
    Note that $\cL[\fs;\ft]+D_1\cH[\fs;\ft]=\overline{\cBs_{\fs}\cC_\ft(R_1)}$. Consider the map $\pi:U_w\to\cS_w^{(\Fq)}$ defined by $\pi(P+D_1Q)=P$, then \Cref{lem:U} shows $\im\pi=\sL_w$ and $\ker\pi=D_1\cR_w$.
\end{proof}

\subsection{Fixed relation}
Now let $V$ be $\cR$ or $U$, and use $V_w^{\text{Fix}}$ to denote the set consisting of $P\in V_w$ such that $[P;0]\in\br_w$ (such $[P;0]$ is called a fixed relation).

\begin{theorem}\label{thm:fix}
    For all $w\in\bbZ_+$,\[\dim_{\Fq}V_w^{\text{Fix}}=\sum_{n=1}^{w}\dim_{\Fq}V_{w-n}.\]
\end{theorem}
\begin{proof}
    For $1\leq n\leq w$, set \[V_{w,n}^{Fix}=\{P\in V_w^{Fix}~|~\rho_k(P)=0\text{ for all }k<n\},\]

    Given a $P_n\in V_{w-n}$, using \Cref{thm:BCP} for $V=\cR$ and \Cref{lem:U} for $V=U$, there exists a binary relation $[A_n;B_n]$ such that $P_n=A_n+B_n$ (Note that such binary relation is not unique, so we should fix it for a basis to ensure well-defined). Then
    \[\cB_{(n)}[A_n;B_n]=\bigl[(n,P_n)+\prescript{n+}{}{B_n};0\bigr],\]
    hence $(n,P_n)+\prescript{n+}{}{B_n}\in V_{w,n}^{Fix}$. This induces an injection $\tau_{w-n}:V_{w-n}\to V_{w,n}^{Fix}/V_{w,n+1}^{Fix}$. 
    
    We show that $\tau$ has an inverse: For a $P=\sum\limits_{k=n}^w(k,P_k)\in V_{w,n}^{Fix}$, $P_n\in V_{w-n}$. By definition, we have \[\sum_{k=n}^wSi_d(k,P_k)=\sum_{k=n}^wl_d^{-k}\Si_{<d}(P_k)=0.\]
    Multiplied by $l_d^{n}$, we have
    \[\Si_{<d}(P_n)+\sum_{k=n+1}^wl_d^{n-k}\Si_{<d}(P_k)=0.\]
    Let $d\to\infty$, note that $l_d^{-1}\to0$, we have $\Li(P_n)=0$, i.e. $P_n\in V_{w-n}$. Then $P\mapsto P_n$ define a linear map from $V_{w,n}^{Fix}\to V_{w-n}$. It factors through $V_{w,n}^{Fix}/V_{w,n+1}^{Fix}$ and is actually the inverse of $\tau_{w-n}$.

    Now we know $\tau_{w-n}$ is an isomorphism. Then we are done by recursion.
\end{proof}

\subsection{Proof of \Cref{thm:dimension}}
\begin{lemma}\label{lem:SRU}
    For all $w\geq q$,
    \[\cS_{w-q}^{(\Fq)}/\cR_{w-q}^{\text{Fix}}\cong U_w/U_w^{\text{Fix}}\]
\end{lemma}
\begin{proof}
    Define the map $\cS_{w-q}^{(\Fq)}\to U_w/U_w^{\text{Fix}}:P\to\overline{\cC_{P}(R_1)}+U_w^{\text{Fix}}$, \Cref{lem:U} and the proof of \Cref{thm:fix} imply it is surjective (In fact, the proof of \Cref{thm:fix} shows that $U_w^{\text{Fix}}$ is generated by these $\overline{\cBs_{\fh}\cC_{\ft}(R_1)}$ with $\fh\neq\emptyset$). We just need to show the kernel of this map is $\cR_{w-q}^{\text{Fix}}$.

    If a $P\in\cS_{w-q}^{(\Fq)}$ satisfies $\overline{\cC_{P}(R_1)}\in U_w^{\text{Fix}}$, then we have
    \[\Si_d(q)\Si_{<d+1}(P)+D_1\Si_{d}\bigl(1,(q-1)\Ss P\bigr)=0.\]
    On the other hand, $\cC_{P}(R_1)$ is a binary relation, we have
    \[\Si_d(q)\Si_{<d+1}(P)+D_1\Si_{d+1}\bigl(1,(q-1)\Ss P\bigr)=0.\]
    So we have $\Si_{d}\bigl(1,(q-1)\Ss P\bigr)=\Si_{d+1}\bigl(1,(q-1)\Ss P\bigr)$ for all $d$. But $\Si_{d}\bigl(1,(q-1)\Ss P\bigr)=0$ for $d$ small enough, so we have $\Si_{d}\bigl(1,(q-1)\Ss P\bigr)=\Si_d(1)\Si_{<d}(q-1)\Si_{<d}(P)=0$ for all $d$. Hence $P\in\cR_{w-q}^{\text{Fix}}$.

    For the converse, it is easy to see if $P\in\cR_{w-q}^{\text{Fix}}$, then $\overline{\cC_{P}(R_1)}\in U_{w}^{\text{Fix}}$.
\end{proof}

Above lemma implies
\[\dim_{\Fq}\cR_{w-q}^{\text{Fix}}=\dim_{\Fq}\cS_{w-q}^{(\Fq)}-\dim_{\Fq}U_w+\dim_{\Fq}U_w^{\text{Fix}}.\]
By \Cref{cor:dimU} ,\Cref{thm:fix} and \Cref{lem:decomposition}, recall we use $r_w$ to denote $\dim_{\Fq}\cR_w$, we have
\[r_w=\sum_{n=1}^qr_{w-n}+|\cI_{w-q}|-|\cI_w-\BTTw|+\sum_{n=1}^w|\cI_{w-n}-\mathcal{B}_{w-n}^{(\textup{TT})}|.\]

Multiplied by $x^w$, summing over all $w\geq1$, use the following generating functions
\[\sum_{w\geq0}|\cI_w|x^w=1+\sum_{w\geq1}2^{w-1}x^w=\frac{1-x}{1-2x},\]
\[\quad\sum_{w\geq0}|\BTTw|x^w=1+\sum_{w=1}^{\infty} d_w^{(K)} x^w = 1+\frac{x(1-x^{q-1})}{1-2x+x^{q+1}},\]
and the following formulation
\[\sum_{w\geq0}(\sum_{k=0}^wa_k)x^w=\frac{1}{1-x}\sum_{w\geq0}a_wx^w.\]
Complicated but direct calculation implies  \Cref{thm:dimension}

\section{Another proof of completeness}\label{sec:another}
In this section, we give another proof of \Cref{thm:BCP}, which does not depend on the reduction of \Cref{thm:BCPHL} (Hence is relatively independent of \cite{HL2026}). The main ideas are similar to the ones in \cite{IKN2026}. We will admit some results in \Cref{sec:proof}.

\subsection{Main idea of proof}

For a $P\in\cS_w^{(\Fq)}$, we use $\min_{\prec}P$ to denote the smallest tuple in $P$ in the (depth,lex)-order.
\begin{lemma}\label{lem:minprefix}
    For $\fs,\ft,\fh,\fu\in\cI$ with $\fs\neq\emptyset$,
    \begin{enumerate}
        \item $\min_{\prec}\overline{\cBs_{\fh}\cP(\fs)}=(\fh,\prescript{q+}{}{\fs}^{+1},q-1)$. We call $\overline{\cBs_{\fh}\cP(\fs)}$ the relation of TYPE I, and $(\fh,\prescript{q+}{}{\fs}^{+1},q-1)$ the tuple of type i;
        \item For $\ft\neq\emptyset$, $\min_{\prec}\overline{\cBs_{\fh}\cP[\fs;\ft]}=(\fh,\prescript{q+}{}{\fs}^{+1},\ft^{+(q-1)})$. We call $\overline{\cBs_{\fh}\cP[\fs;\ft]}$ the relation of TYPE II, and $(\fh,\prescript{q+}{}{\fs}^{+1},\ft^{+(q-1)})$ the tuple of type ii;
        \item For $\fu\neq\emptyset$, $\min_{\prec}\overline{\cBs_{\fh}\Bigl(\cP\bigl[\fs;(\fu^{+1},\ft)\bigr]-\cP\bigl[(\fs^{+1},\fu);\ft\bigr]\Bigr)}=(\fh,\prescript{q+}{}{\fs}^{+1},\fu^{+q},\ft)$. We call \\ $\overline{\cBs_{\fh}\Bigl(\cP\bigl[\fs;(\fu^{+1},\ft)\bigr]-\cP\bigl[(\fs^{+1},\fu);\ft\bigr]\Bigr)}$ the relation of TYPE III, and $(\fh,\prescript{q+}{}{\fs}^{+1},\fu^{+q},\ft)$ the tuple of type iii.
    \end{enumerate}
\end{lemma}
\begin{proof}
    Direct calculation implies desired result.
\end{proof}

By using \Cref{cor:P}, it suffices to show there is no $P\in \cR_w$ such that $\min_{\prec}P$ is not of type i, ii or iii (We call such $P$ illegal). In fact, if we prove it, then for any $P\in\cR_w$, we can repeatedly eliminate the minimal term of $P$ by subtracting some relation in LHS of \Cref{lem:minprefix}. This procedure have to become $0$ and stop after finitely many steps, then we prove that $P$ is a linear combination of BCP relations.

Assume there exists an illegal $P\in \cR_w$. We may choose the $P$ satisfying $\min_{\prec}P$ is maximal among all illegal $P$, then $\ind P$ must be $0$ or legal, because $\min_{\prec}P\prec\min_{\prec}\ind P$ if $\ind P\neq0$. Now we show that $\min_{\prec}P$ is not maximal to derive a contradiction. We will show that $\ind P$ can not be eliminated by subtracting the relations of TYPE I, II or III, then there will remain an illegal relation with minimal tuple greater than $\min_{\prec}P$.

\subsection{Choose tuples}\label{subsec:tuples}

By \Cref{lem:decomposition}, and at most multiplied by a scale, we may write $P=\cL[\fs;\ft]+\sum a_i\cL[\fs_i;\ft_i]$ with $a_i\in\Fq$, $\fs,\fs_i\in[q]^{\bullet}$, and $(\fs,\prescript{q+}{}{\ft})=\min_{\prec}P$ (Let $\prescript{n+}{}{\emptyset}=(n)$ here). Then $\ind P=\cH[\fs;\ft]+\sum a_i\cH[\fs_i;\ft_i]$. Note that \[\min\nolimits_{\prec}\cH[\fs;\ft]=\begin{cases}
    (\fs^{+1},\ft^{+(q-1)}),&\fs\neq\emptyset;\\
    (1,\ft^{+(q-1)}),&\fs=\emptyset.
\end{cases}\]

Now for the case $\fs\neq\emptyset$, set $\dep(\fs,\prescript{q+}{}{\ft})=d$, we choose the tuples in $\ind P$ of the form $(\fs^{+1},\fu)$ with $\dep(\fs^{+1},\fu)=d$ for some $\fu$.

In $\cH[\fs;\ft]$, write $\ft=(t_1,\cdots,t_l)$, such tuples are of form $(\fs^{+1},t_1,\cdots,t_j+q-1,\cdots,t_l)$ for some $j$.

For other $\cH[\fs_i;\ft_i]$, existence of such tuples implies $\dep(\fs_i,\prescript{q+}{}{\ft_i})= d$. In fact, if $\dep(\fs_i)+\dep(\ft_i)>d$, then each tuple of $\cH[\fs_i;\ft_i]$ has depth greater than $d$; $(\fs,\prescript{q+}{}{\ft})\prec(\fs_i,\prescript{q+}{}{\ft_i})$ ensures $\dep(\fs_i)+\max\{1+\dep(\ft_i)\}\geq d$. 

Moreover, we have $\fs_i= \fs$ or $\fs_i=(\fs^{+1},\fu)$ for some $\fu$: If $\dep(\fs_i)<\dep(\fs)$, the tuples of $\cH[\fs_i;\ft_i]$ with depth $d$ have the forms of $(\fs_i^{+1},t_1^{i},\cdots,t_j^{i}+q-1,\cdots,t_m^{i})$, where $\ft_i=(t_1^{i},\cdots,t_m^{i})$. If these tuples have the forms of $(\fs^{+1},\fu)$, let $\fs_i=\bigl(\fs_1^{i},\cdots,\fs_{\dep(\fs_i)}^{i}\bigr)$, $\fs=\bigl(\fs_1,\cdots,\fs_{\dep(\fs)}\bigr)$, then we have $\fs_k^{i}=\fs_k$ for $k<\dep(\fs_i)$ and $\fs_{\dep(\fs_i)}^i+1=\fs_{\dep(\fs_i)}$ by comparing components, which contradicts $(\fs,\prescript{q+}{}{\ft})\prec(\fs_i,\prescript{q+}{}{\ft_i})$. Now we have $\dep(\fs_i)\geq\dep(\fs)$, comparing components derives $\fs=\fs_i$ or $\fs_i=(\fs^{+1},\fu)$.

For the case $\fs=\emptyset$, set $\dep(\prescript{q+}{}{\ft})=d$, we choose the tuples of form $(1,\fu)$ with $\dep(1,\fu)=d+1$. In $\cH[\fs;\ft]$, such tuples are of forms $(1,t_1,\cdots,t_j+q-1,\cdots,t_l)$ for some $j$ with $\ft=(t_1,\cdots,t_l)$.

For other $\cH[\fs_i;\ft_i]$, the same discussion as the above implies that if $\cH[\fs_i;\ft_i]$ has such tuples, then $\fs_i=\emptyset$, $\dep(\ft_i)=d$; or the first component of $\fs_i$ is 1, $\dep(\fs_i)\geq2$, and $\dep(\fs_i,\prescript{q+}{}{\ft_i})=d+1$.

\subsection{Tuples can not be eliminated}
Now we show we can not eliminate the tuples we choose in \Cref{subsec:tuples} to finish our proof. There are two cases:
\begin{enumerate}
    \item [I] $\ft=(1,1,\cdots,1,t)$ or $\emptyset$, i.e. all components of $\ft$ except the last one are $1$ or $\ft$ is empty set.
    \item [II] $\ft\neq\emptyset$ and has components not equal to $1$ before the last one.
\end{enumerate}

\vspace{6pt}
\noindent \textit{\large Case I}: For $\ft\neq\emptyset$, we consider the tuple $(\fs^{+1},\ft^{+(q-1)})=(\fs^{+1},1,1,\cdots,t+q-1)$. It is the minimal tuple of $\cH[\fs;\ft]$. Similar argument as in the previous subsection implies that other $\cH[\fs_i;\ft_i]$ can not generate this tuple. So comparing components implies $(\fs^{+1},\ft^{+(q-1)})=\min_{\prec}\ind P$. However, it is easy to see $(\fs^{+1},\ft^{+(q-1)})$ is not of type i, ii or iii. We are done.

For $\ft=\emptyset$, the same argument works.

\vspace{6pt}
\noindent \textit{\large Case II}: This case will be more complicated.

First, fix a tuple $(v_{1},\cdots,v_n,b)$ with $v_i\leq q$ for all $i$ and $b\leq q-1$. For tuple $(x_1,\cdots,x_{n+1})$ with $x_{n+1}\equiv b\mod{q-1}$, we define $\lambda(x_1,\cdots,x_{n+1})\in\Fq$ by induction: if $1\leq x_i\leq q$ for all $1\leq i\leq n$ and $x_{n+1}=b$, then let \[\lambda(x_1,\cdots,x_{n+1})=\begin{cases}
    1, ~&(x_1,\cdots,x_{n+1})=(v_{1},\cdots,v_n,b);\\
    0, ~&(x_1,\cdots,x_{n+1})\neq(v_{1},\cdots,v_n,b).
\end{cases}\]
We introduce a relation $\leftarrow$. For $\fx=(x_1,\cdots,x_{n+1})$ and $\fy=(y_1,\cdots,y_{n+1})$, let $r$ be the greatest number such that $x_r>q$ if $r\leq n$, or $x_r>b$ if $r=n+1$. We say $\fy\leftarrow\fx$ if there exists an $k<r$ such that $\fy=\bigl(x_1,\cdots,x_k+(q-1),\cdots,x_r-(q-1),\cdots,x_{n+1}\bigr)$. For tuples outside the preceding case, define
\[\lambda(\fx)=\begin{cases}
    0, &\text{if there is no }\fy\text{ such that }\fy<\fx;\\
    -\sum\limits_{\fy\leftarrow\fx}\lambda(\fy),&\text{otherwise}.
\end{cases}\]

Now we set $\fx^{(j)}=(x_1,\cdots,x_j+q-1,\cdots,x_{n+1})$. Then by definition, we have
\[\sum_{j=1}^{n+1}\lambda(\fx^{(j)})=0.\]

For any $i\leq n$ such that $x_i>1$, we have
\[\sum_{j=1}^{i}\lambda(\fx^{(j)})=0.\]
In fact, if there is no $r>i$ such that $x_r>q$ if $r\leq n$, or $x_r>b$ if $r=n+1$, then definition implies the equality. If there exists such $r$, choose the greatest one, let $\fy_k=\bigl(x_1,\cdots,x_k+(q-1),\cdots,x_r-(q-1),\cdots,x_{n+1}\bigr)$. Then
\[\lambda(\fx^{(j)})+\sum_{k=1}^{r-1}\lambda(\fy_k^{(j)})=0,\]
hence
\[\sum_{j=1}^{i}\lambda(\fx^{(j)})=-\sum_{k=1}^{r-1}\sum_{j=1}^i\lambda(\fy_k^{(j)}).\]
By induction on $\sum_{j=1}^{n+1}jx_j$, we can prove it.

Combine these two equalities, we have
\[\sum_{j=i+1}^{n+1}\lambda(\fx^{(j)})=0.\]

We now begin our discussion: In $P=\cL[\fs;\ft]+\sum a_i\cL[\fs_i;\ft_i]$, we may assume we can not eliminate any tuples in $P$ by subtracting relations of TYPE I, II and III (Because subtracting them does not change the minimal tuple of $P$).

For $\ft=(t_1,\cdots,t_l)$, choose $r$ such that $t_j=1$ for all $j\leq r$, and $t_{r+1}\geq2$ (possibly $r=0$). By assumption, $0\leq r \leq l-2$. Now we consider the $\cL[\fs_i;\ft_i]$ (including $\cL[\fs;\ft]$) satisfying $\fs_i=\fs$, $\ft_i=(t_1^{i},\cdots,t_l^{i})$ with $t_j^i=1$ for all $j\leq r$ and $t_{r+1}^i\geq2$. We choose the one whose $t_{r+1}^i$ is maximal among theirs, denote it by $\cL[\fs;\ft']$. We write $\ft'=(\underbrace{1,\cdots,1}_{r~1's},m,v_{1},\cdots,v_n,b)$, and set $\fu=(\underbrace{1,\cdots,1}_{r~1's},m-1)$. Then $\min_{\prec}\cL[\fs;\ft']=(\fs,\prescript{q+}{}{\fu}^{+1},v_1,\cdots,v_n,b)$.

By our assumption that tuples of $P$ can not be eliminated. $(\fs,\prescript{q+}{}{\fu}^{+1},v_1,\cdots,v_n,b)$ has no form of type i, ii or iii. By comparing the components, we have $v_i\leq q$ for $1\leq i\leq n$, $b\leq q-1$; and if $b=q-1$, then $n\geq1$ and $v_n=1$. This argument also works for the previous $\cL[\fs_i;\ft_i]$ we consider.

We consider all the tuples in $\ind P$ of forms
\[(\fs^{+1},\underbrace{1,\cdots,1}_{r~1's},m+q-1,x_1,\cdots,x_{n+1})\quad\text{with }x_{n+1}\equiv b\mod{q-1},\]
where we set $\emptyset^{+1}=1$.

Define an invariant $S=\sum\lambda(x_1,\cdots,x_{n+1})c(x_1,\cdots,x_{n+1})\in\Fq$, where $c(x_1,\cdots,x_{n+1})$ is the coefficient of $(\fs^{+1},\underbrace{1,\cdots,1}_{r~1's},m+q-1,x_1,\cdots,x_{n+1})$ in $\ind P$. Assume the coefficient of $\cL[\fs;\ft']$ in $P$ is $c$, we prove $S=c$:

For the $\cL[\fs_i;\ft_i]$ in $P$ with $\fs_i=\fs$, if $\cH[\fs;\ft_i]$ has the tuples we consider, we have shown that $x_i\leq q$ for all $1\leq i \leq n$ and $x_{n+1}\leq q-1$. By comparing the component, we have $\lambda(x_1,\cdots,x_{n+1})=1\iff\ft_i=\ft'$, and others are $0$. Hence these tuples do not contribute to $S$ except $\cH[\fs;\ft']$ contributing $c$.

For the $\cL[\fs_i;\ft_i]$ in $P$ with $\fs_i\neq\fs$, the similar argument as in previous subsection implies $\fs_i=(\fs^{+1},\underbrace{1,\cdots,1}_{r~1's},m+q-2)$. So the $\cL[\fs_i;\ft_i]$ contribute $\sum\limits_{j=1}^{n+1}\lambda(\ft_i^{(j)})=0$.

Combine above discussion, we get $S=c$.

Finally, we show subtracting relations of TYPE I, II or III does not change $S$. We just need to show it for special relations: We first define $G=\ind P$, if $G$ is illegal, we have done. If not, we subtract a multiple of relation corresponding to the minimal tuple of $G$ to eliminate it, we assume this relation is $\overline{\cBs_{\fh}\cP[\fu;\fv]}$ of TYPE II, other two cases can be treated similarly. Then we have \[\dep(\fh)+\dep(\fu)+\max\{\dep(\fv),1\}=\begin{cases}
    \dep(\min_{\prec}P), &\fs\neq\emptyset;\\
    \dep(\min_{\prec}P)+1, &\fs=\emptyset.
\end{cases}\] At this time, the desired tuples has two possibilities: it can come from
\[\bigl(\fh,\prescript{q+}{}{\fu}^{+1},(q-1)\Ss\fv\bigr),\]
or come from
\[\bigl(\fh^{+1},(q-1)\Ss(\fu,\prescript{q+}{}{\fv})\bigr).\]

If it come from the former. By comparing components with the current minimal tuples, we have \[(\fh,\prescript{q+}{}{\fu}^{+1})=(\fs^{+1},\underbrace{1,\cdots,1}_{r~1's},m+q-1),\] or \[(\fh,\prescript{q+}{}{\fu}^{+1})=(\fs^{+1},\underbrace{1,\cdots,1}_{r~1's},m+q-1,x_1,\cdots,x_u),\] for some $1\leq u\leq n$.

When $\fv\neq\emptyset$, in the former case, $\overline{\cBs_{\fh}\cP[\fu;\fv]}$ contributes
\[\sum_{j=1}^{n+1}\lambda(\fv^{(j)})=0;\]
in the latter case, it contributes
\[\sum_{j=u+1}^{n+1}\lambda\bigl((x_1,\cdots,x_u,\fv)^{(j)}\bigr)=0.\]

When $\fv=\emptyset$, $(q-1)\Ss\emptyset=(q-1)$, hence $b=q-1$. By previous discussion, we have $v_n=1$. But it must be \[(\fh,\prescript{q+}{}{\fu}^{+1})=(\fs^{+1},\underbrace{1,\cdots,1}_{r~1's},m+q-1,x_1,\cdots,x_n),\] as this time. Note that $x_n\geq2>v_n$, and if $x_n>q$, then $x_n-(q-1)\geq2>v_n$. By definition, $\lambda(x_1,\cdots,x_n,q-1)=0$. It also has no contribution to $S$.

If it come from the latter. By comparing components, we have 
\[\fh^{+1}=(\fs^{+1},\underbrace{1,\cdots,1}_{r~1's},m+q-1),\] or \[\fh^{+1}=(\fs^{+1},\underbrace{1,\cdots,1}_{r~1's},m+q-1,x_1,\cdots,x_u),\] for some $1\leq u\leq n$.

The same argument implies zero contribution.

Now if $G-k\overline{\cBs_{\fh}\cP[\fu;\fv]}$ is illegal, where $k$ is the coefficient of $\min_{\prec}G$, we are done. Otherwise use $G-k\overline{\cBs_{\fh}\cP[\fu;\fv]}$ replace $G$, then $\min_{\prec}G\prec\min_{\prec}(G-k\overline{\cBs_{\fh}\cP[\fu;\fv]})$ and repeat above procedure. We have proved $S$ is still equal to $c\neq0$. So $G$ can not become $0$. After finitely many steps, $G$ will become illegal with $\min_{\prec} P\prec \min_{\prec}G$, which derives a contradiction.

\section{Remarks and declarations}

We give some remarks on the proof in this paper to explain the main ideas. Some curious and subtle techniques and constructions can be explained from these ideas, but this is by no means easy.

\Cref{sec:proof} basically follows the ideas in \cite{HL2026}: we want to use the nilpotentency of $\ind$ on $\cR_w$, which can be use to prove certain subspaces of $\cR_w$ are actually the whole space. To do this, establishing \Cref{thm:BCPHL} and studying $\cH(\ker\cL)$ are natural:

Informally, we can think $\ind=\cH\circ\cL^{-1}$, then $\ind^{-1}=\cL\circ\cH^{-1}$. So if we want to find a space $V$ such that $\ind^{-1}(V)\subseteq V$, the simplest idea is making $\cL\circ\cH^{-1}(V)=0$, hence $V=\cH(\ker\cL)$.

However, unlike the arguments in \cite{HL2026} that introduce the so-called $\operatorname{gInd}$ operator and use more linear algebra methods, the space $\BCP_w$ allows us to get desired results by using more techniques about tuples and their operations.

The construction of $\cP(\fs;\ft)$ seems to be mysterious, but directly calculating the LHS of \Cref{thm:CP} basically explains how to find it.

Another key construction in the proof is the RHS of \Cref{thm:KerH}. Unfortunately, the author has no feasible idea for it.

\vspace{6pt}
\Cref{sec:dimension} gives a new method to calculate dimensions: We study the relations among certain spaces instead of detailed linear relations.

The core idea seems to be \Cref{lem:SRU}, which leads to the final dimension formula. However, \Cref{lem:U} may be more essential: We want to calculate $\dim_{\Fq}\BCP_w$, so considering the linear relations among $\cBs$ and $\cC$ is necessary. \Cref{lem:U} ensures that we can obtain some information from $U_w$. Since $\{\overline{\cBs_{\fh}\cC_{\ft}(R_1)}~|~\fh\neq\emptyset\}$ carries the information of $U_w^{Fix}$. So $\overline{\cC_\ft(R_1)}$ will tell us information about $U_w/U_w^{Fix}$, which may be the starting point of \Cref{lem:SRU}. 

\vspace{6pt}
\Cref{sec:another} almost completely follows the path of \cite{IKN2026}. The readers who are familiar with the technique of \cite{IKN2026} may not find it difficult. There are essentially no new ideas but new complicated constructions, such as \Cref{lem:minprefix} (iii) and $\lambda$.

\vspace{6pt}
The author would also like to give some explanations for the BCP relations. The $\cBs$ and $\cC$ operators, due to Todd in \cite{Todd2018} , are ready-made, which can be understood as operations adding prefixes and suffixes to a tuples. If we study the relations computed by computers, after deleting the possible prefixes and suffixes, we may find there should exist other two operations: add $1$ to the second component or insert $1$ before the second component, which are actually the $\cD_0$ and $\cD_1$ operators in \Cref{lem:D}. Todd's idea is that all $K$-relations come from the fundamental relation $R_1$, so we similarly hope that all $\Fq$-relations come from \Cref{eq:basic}, which is the unique $\Fq$-relation with the smallest weight.

\vspace{6pt}
There are some further problems which may be worth discussing:

First, the difference of \Cref{thm:Toddconj} and \Cref{thm:BCP} is what $\cBs\circ\cC$ acts on. It is easy to see $\cP(\fs)=\sum a\cBs_{\fu}\cC_{\fv}(R_1)$ for any $\fs$ and some $a\in K,\fu,\fv\in\cI$. The problem is whether we can obtain more precise description for these $a,\fu,\fv$, which may uncover how to descend from $K$ to $\Fq$.

Another problem is about basis. Can we find a basis with a simple form for $\BCP_w$? It is equivalent to study the linear independence among $\overline{\cBs_{\fu}\cC_{\fv}\cP(\fs)}$ or $\overline{\cBs_{\fu}\cP(\fs;\ft)}$, for which there have been some results when $w<3q$ in \cite{IKN2026}.

The author believes these problems may not be too difficult, but decides to leave them for others’ time or token. The readers with intelligence or artificial intelligence may give it a try.

\vspace{6pt}
\noindent\textbf{About AI use.}
The author provided the construction of $\cD_0,\cD_1$ and some relevant research notes for ChatGPT (GPT-6 Astra, OpenAI), then it generated all mathematical results in this paper.

However, you can still find a few traces of human presence in a flood of AI-generated outputs: The author rewrites and verifies the proofs, and simplifies some of them. The previous remarks are written independently by the author (The author hopes that they will be helpful to the readers, if there are any).

Doubao (ByteDance) was also used to assist the author in writing, especially providing the first sentence of the preceding paragraph.

Since the results were generated by AI, there may be some results not cited properly, and the author apologizes to the relevant authors. Please contact the author for revisions if such cases are found.

\vspace{6pt}
\noindent\textbf{About verification and formalization.}
The author has independently verified all proofs in this paper.

This paper will not be submitted. So if readers find any errors, or find the proofs correct after careful examination, the author would appreciate being notified.

ChatGPT (GPT-6 Astra, OpenAI) also generated the formal proof of this paper in Lean 4, which can be found in \url{https://github.com/jinyuanhu2004/On-the-Fq-linear-relations.git}.

\vspace{6pt}
\noindent\textbf{About credit.}
As for credits, the author has no idea who this belongs to. All he knows is it's definitely not him.

Maybe it can belong to Coco, an extremely cute lady cat, and Bella, a wonderfully beautiful lady Pomeranian.

\begin{acknowledgements}
The author is deeply grateful to ITZY for their encouragement and support. The author also thanks ITZY, TWICE, LE SSERAFIM, NMIXX and MEOVV for staying up late with him while writing this paper.
\end{acknowledgements}

\vspace*{3mm}
\begin{flushright}
\begin{minipage}{148mm}\sc\footnotesize
J.\,Hu: School of Mathematical Sciences, Fudan University, Shanghai, China \\
{\it E-mail address}: \href{mailto:jinyuanhu2004@gmail.com}{{\tt jinyuanhu2004@gmail.com}}
\vspace*{3mm}
\end{minipage}
\end{flushright}

\end{document}